\documentclass[12pt]{amsart}

\usepackage{tasks,comment}
\usepackage{hyperref}
\usepackage{blindtext}
\usepackage{setspace}
\usepackage[utf8]{inputenc}

\usepackage{graphicx}
 
\usepackage{subcaption}
\usepackage[utf8]{inputenc}
\usepackage[export]{adjustbox}
\usepackage{wrapfig}

\usepackage{marginnote}
\usepackage{datetime}
\usepackage{color}
\usepackage{amsrefs}
\usepackage{amsmath}
\usepackage{mathtools}
\usepackage{graphicx}  
\usepackage[top=3.5cm, bottom=2.5cm, outer=2.5cm, inner=2.5cm]{geometry}% margines
\usepackage{hyperref}
\hypersetup{linkcolor=red,
citecolor=black
}
\usepackage{amsmath}
\newtheorem{thm}{Theorem}[section]

\newtheorem{lemma}[thm]{Lemma}
\theoremstyle{definition}

\newtheorem{Alg}[thm]{Algorithm}

\usepackage{hyperref}
\usepackage{hyperref}
\hypersetup{
	colorlinks = true,
	linkcolor = {red},
	citecolor = {blue}
}

\theoremstyle{remark}
\newtheorem{re}[thm]{Remark}

\newcommand{\RNum}[1]{\uppercase\expandafter{\romannumeral #1\relax}}

\numberwithin{equation}{section}

\usepackage[utf8]{inputenc}
\usepackage[english]{babel}
\usepackage{paralist}
\usepackage{amsthm}
\usepackage{marginnote}

\usepackage{todonotes}
 
\theoremstyle{definition}

\theoremstyle{remark}
\newtheorem*{remark}{Remark}

 \usepackage{bm,amsmath}

\newcommand{\RN}[1]{%
  \textup{\uppercase\expandafter{\romannumeral#1}}%
}

\newcommand{\Rey}{\mathcal{R}e }

\newcommand{\bff}{\mathbf{f}}
\newcommand{\bfu}{\mathbf{u}}
\newcommand{\bfv}{\mathbf{v}}
\newcommand{\bfw}{\mathbf{w}}
\newcommand{\bfX}{\mathbf{X}}

\newcommand{\bfV}{\mathbf{V}}

\newcommand{\bfx}{\mathbf{x}}
\newcommand{\bk}{\mathbf{k}}

\newenvironment{Alirev}{\color{blue}}{\color{black}}
\newcommand{\BBB}{\begin{Alirev}}
\newcommand{\FFF}{\end{Alirev}}
\newenvironment{Mikerev}{\color{blue}}{\color{black}}
\newcommand{\MMM}{\begin{Mikerev}}
\newcommand{\JJJ}{\end{Mikerev}}

\allowdisplaybreaks[4]
\usepackage[utf8]{inputenc}

\title[IMEX-FEM schemes for nonlinear  data
assimilation]{Finite element nonlinear nudging for the Navier-Stokes equations}

\author[Carlson]{Elizabeth Carlson}
\author[Fang]{Rui Fang$^*$}
\author[Pakzad]{Ali Pakzad}

\thanks{$^*$Corresponding author.}

\address[E. Carlson]{Department of Mathematics, Oregon State University}
\email{carleliz@oregonstate.edu}

\address[R. Fang]{Department of Mathematics, The Ohio State University}
\email{fang.1211@osu.edu}

\address[A. Pakzad]{Department of Mathematics, California State University Northridge}
\email{pakzad@csun.edu}
\date{\today}


\begin{document}
\begin{abstract}
Nonlinear nudging augments the classical linear nudging (continuous data assimilation) term with a nonlinear feedback that strengthens as the observed error decreases, and yields super‑exponential convergence to the true state. In this paper, we study nonlinear nudging for the $3$-dimensional  incompressible Navier–Stokes equations with the finite element method, which allows complex geometries and physical boundary conditions. We propose an implicit–explicit (IMEX) finite element scheme with inf‑sup stable elements, prove its unconditional long‑time stability, and show for the spatially semi‑discrete scheme that the error reaches a prescribed tolerance in finite time. Our analysis also suggests that the nonlinear nudging term requires the mesh resolution to be fine enough ($\chi =\mathcal{O}(H^{-1})$, where $\chi$ is the nonlinear nudging intensity and $H$ is the observation mesh size). Further, our error estimates suggest that a larger nonlinear exponent $\gamma$ requires finer observational resolution, which is confirmed by numerical tests. Numerical experiments, including two‑dimensional forced viscous flow in a periodic box, flow past a cylinder, and flow over a flat plate, show that nonlinear nudging synchronizes with the true solution considerably faster than linear nudging. We also show that observations in a subregion suffice to synchronize the whole flow. 
\end{abstract}
\maketitle

\section{Introduction}
In computational fluid dynamics (CFD),  accurate prediction of the future state of an incompressible viscous fluid in a bounded physical domain depends not only on the fidelity of the governing equations but also on the quality of available data. Data assimilation \cite{K03} provides a systematic framework for incorporating discrete observational measurements into numerical simulations in order to reduce modeling and initialization errors \cite{Charney-Use1969}. The objective is to combine sparse \cite{Cibik-Fang-Layton2026}, local \cites{FangPakzad2026,Biswas-Data2021},  noisy \cite{Bessaih-Continuous2015} and possibly time-delayed and discrete \cite{Biswas-Downscaling2019} data with approximate numerical discretizations of the flow equations to improve stability, improve long-term accuracy, and produce more reliable predictions of the flow state and its statistical behavior \cite{K03}.

Current best practice in data assimilation often employs Kalman
filters and their variants \cite{Kalman1960}. Recently, modular nudging algorithms inspired by Kalman filters
have been developed and analyzed for data assimilation
\cite{Cibik-Fang-Layton2026}. Since the 2014 work of Azouani, Olson, and Titi \cite{AOT14}, nudging, also known as
continuous data assimilation (CDA), has developed a strong
mathematical foundation (see, e.g.,
\cites{Biswas_3D, CAO2022103659, leotest}) and has been
implemented in a variety of settings (see, e.g.,
\cites{Pawar-Long2020, CIBIK2025117526, Leo2022, ALT16,
Newey-Whitehead-Carlson2025}).

We begin by considering the dynamics of incompressible viscous flows widely modeled by the incompressible Navier--Stokes equations, and indeed also serve as a fundamental framework for fluid flow prediction:
\begin{equation}
 \label{NSE}
 \begin{split}
\bfu_t+   (\bfu \cdot \nabla) \bfu   - \nu \Delta \bfu +  \nabla p  & = \bff,\\ 
\nabla \cdot \bfu  &= 0, 
\end{split}
\end{equation}
in the physical domain $\Omega \subset \mathbb{R}^2$ or $\mathbb{R}^3$ \cite{FMRT01} with an  initial condition $\bfu_0(x)$. Here,  the true velocity is denoted by $\bfu(\bfx,t)$, and its observed value $I_H(\bfu)$ represents an interpolation of the true state onto a spatially discrete set of measurement points, satisfying \eqref{Ineq:Intrpolate}.   The parameter $H$ denotes the characteristic length scale associated with the observational mesh and measures the spacing between observation points.  In \eqref{NSE}, $p$ denotes the pressure, $\bff$ is the external forcing, and $\nu>0$ is the kinematic viscosity, all of which determine the dynamics of the fluid.

In practice, the initial state $\bfu_0(\bfx)$ may be missing or only partially known. This uncertainty in the initial condition can lead to exponential growth of the solution error in a nonlinear system. To address this difficulty, data assimilation methods directly incorporate observational measurements into the mathematical model. 
Data assimilation falls into two camps: deterministic (nudging-like algorithms, e.g. \cites{Hoke-TheInitialization1976, Stauffer_Seaman_1990, AOT14}) and stochastic (Kalman filter and variants, e.g., \cites{Kalman1960, Evensen-Data2009, hunt2007efficient}, 4DVAR, e.g., \cites{ledimet1986variational, lewis1985use} and particle filter \& variants, e.g., \cites{Gordon_etal_1993, Kitagawa1996}.  In this paper, we focus on a deterministic algorithm known as CDA algorithm \cite{AOT14} for its theoretically predictive power.  It has since been successfully applied to a wide class of dissipative  partial differential equations (PDEs) \cites{Bessaih-Continuous2015, Cao-Algebraic2021, Farhat-Data2020, Jolly-Continuous2019, Pei-Continuous2019, Newey-Whitehead-Carlson2025, CAO2022103659} and provides a computationally efficient alternative to more elaborate statistical data assimilation techniques.  Moreover, in recent years it has been developed and expanded in multiple directions, including but not limited to reducing model errors arising from approximate models (a very incomplete but a list from which additional references can be explored: \cites {Larios_Pei2020,chen2023data,hawkins2023removing,cibik2025data,CibikFang2026Compressible,CibikFang2026Numerical, CIBIK2025117526, Carlson_Hudson_Larios2020, Pachev_Whitehead_McQuarrie_2021concurrent, Newey-Whitehead-Carlson2025}).  It has also, as is common practice in data assimilation, led to new data assimilation algorithms using inspiration from or the hybridization of existing ideas (see, e.g., \cites{Cibik-Fang-Layton2026, Pawar-Long2020, Wang_Carlson_Hoffmann_2026, CLT23}). An algebraic nudging operator was proposed in~\cite{LeoAlgebricDA} that adds continuous data assimilation to existing finite element and legacy codes through a simple diagonal matrix, while retaining long-time stability and synchronization. For a complete review of all extensions of the AOT theory, see citations of \cite{AOT14}.  %In this particular project, we focus on extending the work in \cite{CLT23} in the case of a fully discrete approximation of the Navier-Stokes equations on a physical domain.

Recently, Carlson, Larios, and Titi  \cite{CLT23} introduced a nonlinear modification of the linear nudging algorithm in the form 
\begin{equation}
 \label{NLDANSE}
 \begin{split}
\bfv_t+   (\bfv \cdot \nabla) \bfv   - \nu \Delta \bfv +  \nabla q   = \bff + \chi \,  \mathcal{N} (I_H(\bfu) -  I_H(\bfv)) +   \mu\, I_H(\bfu - \bfv ) &,\\ 
\nabla \cdot \bfv  = 0& ,
\end{split}
\end{equation}
where   $\chi , \mu>0$ are  relaxation parameters and  for $\gamma \in [0,1)$, the nonlinear operator takes the form of 
\begin{equation}
\label{NLterm}
    \mathcal N(\phi) =
\begin{cases}
0, & \|\phi\|_{L^2(\Omega)}=0,\\
\phi\,\|\phi\|_{L^2(\Omega)}^{-\gamma}, & \|\phi\|_{L^2(\Omega)}>0.
\end{cases}
\end{equation}

Under suitable assumptions on the parameters, this nonlinear-nudging algorithm  yields global well-posedness and  was shown to exhibit exponential convergence, and, once the error falls below a certain threshold, a \emph{double-exponential} (super-exponential) convergence rate up to a precision determined by the resolution of the observed data. Computational studies for the one-dimensional Kuramoto--Sivashinsky equation \cite{Larios_Pei2020} and the two-dimensional incompressible Navier--Stokes equations \cite{CLT23} demonstrated that the nonlinear feedback accelerates synchronization of the large scales significantly compared to the linear scheme. The key advantage of nonlinear nudging lies in its ability to enhance transient convergence while retaining long-time stability properties.

All previous simulations of nonlinear nudging were carried out using pseudo-spectral methods on simple periodic domains with uniform meshes. However, many realistic applications involve complex geometries, irregular domains, and non-periodic boundary conditions, where spectral methods are not sufficient.  In such settings, the finite element method (FEM) provides a natural and flexible framework for spatial discretization. While PDE-level convergence theory is essential, discretization introduces additional numerical considerations, including mesh resolution, stability constraints, element choice, and boundary treatment, that may critically affect the success of the data assimilation algorithm. Therefore, it is both theoretically and practically important to investigate the nonlinear-nudging algorithm within a FEM framework, ensuring stability, accuracy, and robustness in geometrically complex and physically realistic configurations.  %Moreover, the original theory had strong, and somewhat unphysical, restrictions on the types of data interpolants one could use.  We get around these restrictions, at the cost of reducing to exponential convergence.

\subsection*{Results of this paper} In this paper, we use inf-sup stable mixed FEMs to approximate \eqref{NLDANSE}. The analysis holds in both two and three dimensions, while the numerical experiments are carried out in two dimensions due to the high computational cost of resolved three-dimensional simulations.

In Section \ref{S3}, we first establish stability of the fully
discrete scheme with first-order implicit–explicit (IMEX) time stepping. The analysis is based on a condition on $H$ and $\mu$: $\mu$ is chosen to satisfy \cite{CIBIK2025117526}
\begin{equation}
\begin{gathered}
  \mu - 2187 C_2^4  \nu^{-3}\|\nabla \bfv^h\|_{L^2(\Omega)}^4 \geq \alpha \mu \geq 0, \\
  \frac{\nu}{4}-\mu C^2_1 H^2 \geq 0.
  \end{gathered}
\end{equation}

In Section~\ref{error-estimates},
we derive error estimates for the semi-discrete scheme and show
that the error reaches a prescribed tolerance in finite time.
When the error is large, the linear nudging term dominates
and drives the error reduction. As the error becomes small,
the nonlinear nudging term becomes dominant and accelerates
the decay.

Finally, in Section~\ref{numerical_test}, we demonstrate the accuracy, efficiency, and robustness of the nonlinear nudging method through numerical experiments. We first verify that the scheme with the second-order backward differentiation formula (BDF2) achieves second-order accuracy in time. We then consider three flow problems: a two-dimensional forced flow in a periodic box, and two flows past an obstacle in a channel, one with a smooth circular obstacle (a cylinder) and one with a thin flat plate, which is a more singular obstacle.
In all tests, nonlinear nudging reduces the velocity error faster than linear nudging. We also study the effect of the exponent $\gamma$, which controls the strength of the nonlinear feedback. A larger $\gamma$ accelerates the early convergence but can lead to a higher error plateau, which is consistent with our error estimates: a stronger nonlinear effect requires finer observational data. Finally, motivated by  \cite{Biswas-Data2021} and \cite{FangPakzad2026},  we extend nonlinear nudging to the case where observations are available only in a subregion of the domain, and show that it remains effective when the observed region contains the obstacle. For the flow past a cylinder, the assimilated solutions also recover the drag and lift on the obstacle, with errors that follow the velocity errors.

\section{Notations and Preliminaries}

We consider a bounded open domain $\Omega \subset \mathbb{R}^d$, where $d=2$ or $3$.  The $L^2(\Omega)$ norm and inner product will be denoted by $\|\cdot\|$ and $(\cdot,\cdot)$, respectively, and $H^1(\Omega)$ seminorm is denoted by $|\cdot|$, while all other norms will be labeled with subscripts. The notation for the FEMs used in CFD is standard and well documented in the literature; see, for example, \cite{layton2008introduction} and \cite{John2016}. For completeness, we briefly review the notation required throughout this manuscript. 

Denote the natural function spaces for velocity 
and pressure, respectively, by
$$
\bfX := (H_0^1 (\Omega))^d =  \{\bfu \in (H^1(\Omega))^d: \nabla \bfu \in (L^2(\Omega))^{d\times d}, \bfu = 0  \text{ on }  \partial \Omega \}
$$
and 
$${Q}:= L_0^2(\Omega) = \{q \in L^2(\Omega): \hspace{3pt}\int_{\Omega} q\  dx= 0\}.$$
The closed subspace of divergence free functions is given by
$$ \bfV :=   \{\bfu \in \bfX: \, (\nabla \cdot \bfu , q ) = 0, \, \, \forall \ q \in Q \},$$
and the dual norm of $\bfX$ will be denoted by $\|\cdot\|_{-1}$. In $\bfX$, we have the Poincar\'e inequality: there exists a constant $C_P$ depending only on $\Omega$ such that 
for any $\phi \in \bfX$,
\[
\|\phi\| \leq C_P \|\nabla \phi\|.
\]

We denote the (explicitly skew symmetrized)  trilinear form $b : \bfX \times \bfX  \times \bfX  \rightarrow \mathbb{R}$ as 
\begin{equation}
b (\bfu, \bfv, \bfw) := \frac{1}{2} (\bfu \cdot \nabla \bfv, \bfw) - \frac{1}{2} (\bfu \cdot \nabla \bfw, \bfv).
\end{equation}
Further, we have
\begin{equation}
b (\bfu, \bfv, \bfw) = (\bfu \cdot \nabla \bfv, \bfw) + \frac{1}{2} ((\nabla \cdot \bfu) \bfv, \bfw).
\end{equation}
We define
\begin{equation}
    K:= \sup_{\bfu,\bfv,\bfw\in \bfV} \frac{b(\bfu,\bfv,\bfw)}{\|\nabla \bfu\|\|\nabla \bfv\| \|\nabla \bfw\|}<\infty.
\end{equation}
\begin{lemma} \label{nonlinear_bound_John}
(\cite{John2016} p. 312) There is a $C_2<\infty$ with
\begin{equation}
\begin{split}
 b(\bfu, \bfv, \bfw)
&\leq C_2 \sqrt{\|
\bfu\|}\sqrt{\|\nabla \bfu\|} \|\nabla \bfv\|\|\nabla \bfw\|,
 \end{split}
\end{equation}
for $\bfu,\bfv,\bfw \in (H^1_0(\Omega))^d$, $d=2$ or $3$.
\end{lemma}

\begin{lemma} (H\"older's and Young's inequalities) For any $\sigma>0$, $1\leq q <\infty$, and $\frac{1}{p}+\frac{1}{q}=1$, the H\"older's and Young's inequalities:
\begin{equation}
    (\bfu,\bfv) \leq \|\bfu\|_{L^p} \|\bfv\|_{L^q}, \ (\bfu,\bfv) \leq \frac{\sigma}{p} \|\bfu\|^p_{L^p} + \frac{\sigma^{-q/p}}{q} \|\bfv\|^q_{L^q}.
\end{equation}

\end{lemma}

\subsection{Discretization preliminaries} Let $\Omega_H$ denote a coarse finite element mesh representing the spatial locations where observational data is available in practice.  To incorporate observational data, we introduce an interpolation operator
\[
I_H:(H^1(\Omega))^d\to (L^2(\Omega))^d,
\]
which maps the true solution onto a coarse observational grid with characteristic mesh size \(H\). We assume that \(I_H\) satisfies the standard approximation and stability properties
\begin{equation}\label{Ineq:Intrpolate}
\begin{split}
\|\phi-I_H\phi\| &\leq C_1H\|\nabla\phi\|,\\
\|I_H\phi\| &\leq C_I\|\phi\|,
\end{split}
\end{equation}
for all \(\phi\in(H^1(\Omega))^d\), where the constants \(C_1\) and \(C_I\) are independent of \(\phi\) and \(H\); see, for example, \cite{BR08}. Examples of interpolation operators satisfying \eqref{Ineq:Intrpolate} include projections onto Fourier modes with wave numbers \(|k|\leq H^{-1}\), local volume-element interpolants, and piecewise constant finite-element interpolants \cites{Biswas-Data2021,AOT14}. In the finite-element setting, \(I_H\) may also be chosen as a piecewise polynomial interpolation operator, such as a piecewise linear interpolant. Numerical evidence further suggests that higher-order interpolation may accelerate synchronization under nudging \cite{JP23}.

In the analysis below, we will consider the special case in which \(I_H\) is the \(L^2\)-orthogonal projection
\begin{equation}\label{l2-project}
(I_H v,w)=(I_H v,I_H w)
\qquad\text{for all }v,w\in(L^2(\Omega))^d.
\end{equation}

The coarse mesh $\Omega_H$ is successively refined to obtain a finer mesh $\Omega_h$, with $h \ll H$. The fine mesh $\Omega_h$ corresponds to the discretization used for the Direct Numerical Simulation (DNS). In particular, $\Omega_h$ is taken to be a conforming triangulation of the domain $\Omega$ with maximum element diameter $h$.  

Let $\bfX^h \subset \bfX$ and $Q^h \subset Q$ denote finite-dimensional velocity and pressure spaces associated with this mesh. We assume that the pair $(\bfX^h, Q^h)$ is inf-sup stable, i.e., it satisfies the discrete Ladyzhenskaya--Babu\v{s}ka--Brezzi (LBB) condition
\begin{equation}\label{LBBH}
\inf_{q^h \in Q^h} \sup_{\bfv^h \in \bfX^h} 
\frac{(q^h, \nabla \cdot \bfv^h)}{\|\nabla \bfv^h\| \, \|q^h\|} 
\geq \beta^h > 0,
\end{equation}
where the constant $\beta^h$ is independent of $h$ as $h \to 0$ (see, e.g., \cite{John2016}). Equivalently, the inf-sup condition can be written as
\[
\beta^h \|q^h\| \leq \sup_{\bfv^h \in \bfX^h} 
\frac{(q^h, \nabla \cdot \bfv^h)}{\|\nabla \bfv^h\|}.
\]
This condition guarantees stability of the discrete velocity and pressure approximations to the Navier--Stokes equations, in the sense that
\[
\|\nabla \bfu^h\|^2 + \|p^h\|^2 
\leq C_h\left( \|\nabla \bfu\|^2 + \|p\|^2 \right).
\]

In finite element discretizations of the Navier--Stokes equations, the incompressibility condition $\nabla \cdot \bfu^h = 0$ is typically enforced in a weak sense rather than pointwise. As a result, the discrete velocity $\bfu^h \in \bfX^h$ is constrained to lie in the discretely divergence-free subspace $\bfV^h \subset \bfX^h$, defined by
\[
\bfV^h := \left\{ \bfv^h \in \bfX^h : (r^h, \nabla \cdot \bfv^h) = 0 \quad \forall \, r^h \in Q^h \right\}.
\]

Motivated by the error equation (1.29) in \cite{girault1979finite}  p. 19 and the discrete inf-sup condition (\ref{LBBH}), also see \cite{cibik2025data}, we assume $\bfV^h$ has approximation properties in (\ref{prop}) similar to $\bfX^h$, for $\bfu\in (H^{m+1}(\Omega))^d\ \cap\ (H_0^1(\Omega))^d$ and $p\in H^{m}(\Omega)$,
\begin{equation}\label{prop}
\begin{aligned}
\inf_{\bfv^h\in \bfV^h}\{\|\bfu-\bfv^h\|+h\|\nabla(\bfu-\bfv^h)\|\}&\leq C_4h^{m+1}|\bfu|_{m+1},
\\
\inf_{q^h\in Q^h}\|p-q^h\|&\leq C_4h^{m}|p|_{m}.
\end{aligned}
\end{equation}

 \begin{lemma}\label{Lem;Seq} Suppose the sequences $\{a_n\}_{n=0}^\infty$ of real number satisfies 
 $$ r a_{n+1} \leq a_n +K. $$
 Let $r > 1$ and $K\geq 0$, then we have 
 $$a_{n+1} \leq a_0\, \left(\frac{1}{r}\right)^{n+1} + \frac{K}{r-1}.$$
 \end{lemma}

\section{Finite element discretization and a first order IMEX-FEM scheme}\label{S3}

We consider a fully discrete finite-element approximation of the nonlinear data assimilation system \eqref{NLDANSE}. The weak form of the nonlinear nudging scheme $(\ref{NLDANSE})$:
\begin{equation}
 \label{weakfNLDANSE}
 \begin{split}
(\bfv_t, \Theta)+   ((\bfv \cdot \nabla) \bfv ,\Theta)  +\nu (\nabla \bfv, \nabla \Theta) - (q, \nabla\cdot \Theta) =(\bff,\Theta) & \\
 + \chi \,  (\mathcal{N} (I_H(\bfu) -  I_H(\bfv)),\Theta) +   \mu\, (I_H(\bfu - \bfv ),\Theta)&,
 \end{split}
\end{equation}
for all $\Theta \in \bfV$. 

We first consider the spatially
semi-discrete approximation of \eqref{weakfNLDANSE}, in which time
remains continuous. Given $\bfv^h(0)=\bfv_0^h\in\bfV^h$, seek
$(\bfv^h,q^h)\in\bfX^h\times Q^h$ such that, for all
$(\Theta^h,r^h)\in\bfX^h\times Q^h$, the following equations hold:
\begin{equation}
 \label{semifNLDANSE}
 \begin{split}
(\bfv^h_t, \Theta^h)+  b(\bfv^h,\bfv^{h},\Theta^{h}) +\nu (\nabla \bfv^h, \nabla \Theta^h) - (q^h, \nabla \cdot \Theta^h) = (\bff,\Theta^h) & \\
 + \chi \,  (\mathcal{N} (I_H(\bfu) -  I_H(\bfv^h)),\Theta^h) +   \mu\, (I_H(\bfu - \bfv^h ),\Theta^h)&.
 \end{split}
\end{equation}

Time integration is performed using a first-order backward Euler scheme. As described in Algorithm \ref{Algorithm}, the convective term is treated semi-implicitly using an IMEX approach: the convecting velocity is evaluated at the previous time level, while the transported velocity is evaluated at the current time level. Thus,
\begin{equation}  (\bfv^h_{n+1}\cdot\nabla)\bfv^h_{n+1}
    \leadsto
    (\bfv^h_n\cdot\nabla)\bfv^h_{n+1}.
\end{equation}
This treatment avoids the need to solve a nonlinear system arising from the convective term at each time step. In contrast, the nonlinear nudging term is retained implicitly at the current time level:
\begin{equation}
    \mathcal{N}_{n+1}
    =
    \frac{I_H(\bfu_{n+1})-I_H(\bfv^h_{n+1})}
    {\|I_H(\bfu_{n+1})-I_H(\bfv^h_{n+1})\|^{\gamma}},
\end{equation}
where $\mathcal{N}_{n+1}=0$ whenever the observational mismatch vanishes.

\begin{Alg}\label{Algorithm}
  Given body force $\bff \in L^{\infty} ((0, \infty), L^2(\Omega))$,  initial condition $\bfv_0^h \in \bfV^h$ , reference solution $\bfu_{n+1}$, and  $(\bfv^h_n, q_n^h) \in \bfX^h \times Q^h$,  compute  $(\bfv_{n+1}^h, q_{n+1}^h) \in \bfX^h \times Q^h$ satisfying
\begin{multline}\label{IMEX-FEM}
\left(\frac{\bfv_{n+1}^{h} - \bfv_{n}^{h}}{\Delta t},\Theta ^{h}\right) + b(\bfv^h_n,\bfv_{n+1}^{h},\Theta^{h}) + \nu \, (\nabla \bfv_{n+1}^{h},\nabla \Theta^{h})  - (q_{n+1}^{h}, \nabla \cdot \Theta^{h})
\\ =  (\bff_{n+1},\Theta^{h})  + \mu \, (I_H(\bfu_{n+1} - \bfv^h_{n+1} ),  \Theta^h)  + \chi \,  \left( \frac{I_H(\bfu_{n+1}) -  I_H(\bfv^h_{n+1})}{\|(I_H(\bfu_{n+1}) -  I_H(\bfv^h_{n+1}))\|^{\gamma}},  \Theta^h\right)
\end{multline}
\begin{equation}\label{incomp}
(\nabla \cdot \bfv_{n+1}^{h},r^{h}) = 0,
\end{equation}
for all $(\Theta^{h},  r^{h}) \in  \bfX^h \times Q^h.$
\end{Alg}
Starting from an arbitrary initial condition $\bfv^h_0$, the approximate solution \(\bfv_{n+1}^h\) is generated using Algorithm \ref{Algorithm}. At each time level \(t_{n+1}\), the quantity \(I_H(\bfu_{n+1})\) denotes the coarse spatial observations of the true state. Because these measurements are assumed to be available when advancing the approximation to \(t_{n+1}\),  then \(I_H(\bfu_{n+1})\) serves as the latest observational input to the algorithm.

For the spatial discretization, we employ the FEM on general meshes, providing the flexibility needed to accommodate complex geometries and boundary conditions. The velocity and pressure spaces are chosen as an inf-sup stable pair \((\bfX^h,Q^h)\) satisfying \eqref{LBBH}. For example, one may use the Taylor–Hood pair \((P_k,P_{k-1})\), with \(k\geq2\). Other LBB-stable finite-element pairs may also be used without substantial changes to the algorithm (see \cites{layton2008introduction, John2016}).

\begin{re}
For computational efficiency, the numerical experiments use a variant of Algorithm \ref{Algorithm} in which the denominator of the nonlinear nudging term
is evaluated at the previous time level, while its numerator remains
at the current time level. Together with the semi-implicit treatment
of convection, this modification yields a linear system at each time step.
\begin{equation}
\begin{gathered}
   \mathcal{N}_{n+1}=  \frac{I_H(\bfu_{n+1}) -  I_H(\bfv^h_{n+1})}{\|(I_H(\bfu_{n+1}) -  I_H(\bfv^h_{n+1}))\|^{\gamma}}  \leadsto \mathcal{N}_{n+1}^*= \frac{I_H(\bfu_{n+1}) -  I_H(\bfv^h_{n+1})}{\|(I_H(\bfu_{n}) -  I_H(\bfv^h_{n}))\|^{\gamma}}.
\end{gathered}
\end{equation} 
\end{re}

We now establish stability and error estimates for Algorithm \ref{Algorithm}, with the nonlinear nudging term treated implicitly. We begin by proving in Theorem \ref{stability} that the algorithm is well-posed globally in time.

\begin{thm}\label{stability}
Let $0\le \gamma<1$, and assume that $I_H$ to be the $L^2$-orthogonal projection
defined in \eqref{l2-project}, 
$U:=\|\bfu\|_{L^{\infty}(0,\infty;L^2(\Omega))}<\infty$,
$F:=\|\bff\|^2_{L^{\infty}(0,\infty;H^{-1}(\Omega))}<\infty$.
Then Algorithm \ref{Algorithm} is unconditionally stable: for all $n\geq 0$,
\begin{equation}\label{stability-final-result}
\|\bfv^h_{n+1}\|^2 \leq
\Big(1+\Delta t\frac{\nu}{C_{FP}^2}\Big)^{-(n+1)}\|\bfv^h_0\|^2
+\frac{M}{\frac{\nu}{C_{FP}^2}},
\end{equation}
where
$M:=\frac{2F}{\nu} + \mu \, U^2 +  \frac{4\chi^2}{\mu}  U^{2-2\, \gamma}  +  \chi \gamma \left(\frac{4\chi (2-\gamma)}{\mu}\right)^{(2-\gamma)/\gamma} $
is independent of $\Delta t$, $h$, and $n$.
\end{thm}

\begin{proof}
We consider the case $0<\gamma<1$. When $\gamma=0$,
Algorithm~\ref{Algorithm} reduces to the standard linear-nudging
scheme. We set $\Theta^{h} = \mathbf{v}^h_{n+1}$ in \ref{Algorithm}. Through the application of the polarization identity, and by the skew-symmetry of the nonlinear term, we reach the following expression:
\begin{equation}\label{Disc1}
\begin{gathered}
\frac{1}{2\Delta t } \left( \| \bfv^h_{n+1}\|^2 -\| \bfv^h_{n}\|^2  \right)  + \nu\, \|\nabla \bfv^h_{n+1}\|^2  + \mu \,  (I_H(\bfv^h_{n+1}), \bfv^h_{n+1}) \\\leq (\bff_{n+1} , \bfv^h_{n+1}) 
 + \mu \, (I_H(\bfu_{n+1}),  \bfv^h_{n+1}) \\
 + \chi \,  \left( \frac{I_H(\bfu_{n+1}) -  I_H(\bfv^h_{n+1})}{\|(I_H(\bfu_{n+1}) -  I_H(\bfv^h_{n+1}))\|^{\gamma}},  \bfv^h_{n+1}\right).
\end{gathered}
\end{equation}
Since $I_H$ is $L^2$-orthogonal projection, 
\begin{equation}
    \mu \,  (I_H(\bfv^h_{n+1}), \bfv^h_{n+1})=  \mu \,\|I_H(\bfv^h_{n+1})\|^2.
\end{equation}
Similarly, by the projection property, and applying Cauchy-Schwarz inequality and Young's inequality, we have
\begin{equation}
    \mu \, (I_H(\bfu_{n+1}),  \bfv^h_{n+1})\leq \frac{1}{2}\mu \|I_H(\bfu_{n+1})\|^2+\frac{1}{2}\mu \|I_H(\bfv^h_{n+1})\|^2.
\end{equation}
For the nonlinear nudging, since $I_H$ is $L^2$-orthogonal projection, Cauchy-Schwarz inequality, sub-additivity, and Young's inequality yield 
\begin{equation}
\begin{gathered}
\chi \,  \left( \frac{I_H(\bfu_{n+1}) -  I_H(\bfv^h_{n+1})}{\|I_H(\bfu_{n+1}) -  I_H(\bfv^h_{n+1})\|^{\gamma}},  \bfv^h_{n+1}\right)   \leq \chi  \|I_H(\bfu_{n+1}) -  I_H(\bfv^h_{n+1})\|^{1-\gamma}\, \| I_H(\bfv^h_{n+1})\|\\
 \leq \chi\, \|I_H(\bfu_{n+1})\|^{1-\gamma}\|I_H(\bfv^h_{n+1})\| + \chi\,    \|I_H(\bfv^h_{n+1})\|^{2-\gamma}\\
 \leq   \frac{2\chi^2}{\mu}  \|I_H(\bfu_{n+1}) \|^{2-2\, \gamma}  +   \frac{ \chi \gamma}{2} \left(\frac{4\chi (2-\gamma)}{\mu}\right)^{(2-\gamma)/\gamma} + \frac{\mu}{4} \| I_H(\bfv^h_{n+1})\|^2.
\end{gathered}
\end{equation}
The first term on the right-hand side is bounded utilizing the  Young's inequality,
\begin{equation}(\bff_{n+1} , \bfv^h_{n+1}) \leq \frac{1}{2\nu}\|\bff_{n+1}\|^2_{-1} + \frac{\nu}{2} \| \nabla \bfv^h_{n+1}\|^2.
\end{equation}
Combine the above bounds for the right hand side of \eqref{Disc1} to arrive at
\begin{equation*}
\begin{split}
\left( \| \bfv^h_{n+1}\|^2 -\| \bfv^h_{n}\|^2  \right)  + \Delta t \, \nu\, \|\nabla \bfv^h_{n+1}\|^2  + \Delta t \, \frac{\mu}{2}\,  \|I_H(\bfv^h_{n+1})\|^2\\
\leq  \Delta t \left( \frac{1}{\nu} \|\bff_{n+1}\|^2_{-1}+\mu\|I_H(\bfu_{n+1})\|^2+ \frac{4\chi^2}{\mu}  \|I_H(\bfu_{n+1}) \|^{2-2\, \gamma}  +  \chi \gamma \left(\frac{4\chi (2-\gamma)}{\mu}\right)^{(2-\gamma)/\gamma}\right).
\end{split}
\end{equation*}
Thus, we have
\begin{equation}
\begin{gathered}
\left( \| \bfv^h_{n+1}\|^2 -\| \bfv^h_{n}\|^2  \right)  + \Delta t \, \nu\, \|\nabla \bfv^h_{n+1}\|^2  + \Delta t \, \frac{\mu }{2}\,  \|I_H(\bfv^h_{n+1})\|^2 \\
\leq  \Delta t \left( \frac{2F}{\nu} + \mu \, U^2 +  \frac{4\chi^2}{\mu}  U^{2-2\, \gamma}  +  \chi \gamma \left(\frac{4\chi (2-\gamma)}{\mu}\right)^{(2-\gamma)/\gamma} \right).
\end{gathered}
\end{equation}
Applying the Poincar\'e inequality to the viscous term and collecting terms gives us
\begin{equation}
\begin{gathered}
\left( 1 + \Delta t \,  \frac{\nu}{(C_{PF})^2} \right)  \| \bfv^h_{n+1}\|^2  \leq   \| \bfv^h_{n}\|^2  \\
+\Delta t \left( \frac{2F}{\nu} + \mu \, U^2 +  \frac{4\chi^2}{\mu}  U^{2-2\, \gamma}  +  \chi \gamma \left(\frac{4\chi (2-\gamma)}{\mu}\right)^{(2-\gamma)/\gamma} \right).
\end{gathered}
\end{equation}
Using Lemma \ref{Lem;Seq}, we have
\begin{equation}
\begin{gathered}
\| \bfv^h_{n+1}\|^2    \leq   \frac{1}{\left( 1 + \Delta t \,  \frac{\nu}{(C_{PF})^2} \right)^{1+n}}\, \| \bfv^h_{0}\|^2\\
+ \frac{ \frac{2F}{\nu} + \mu \, U^2 +  \frac{4\chi^2}{\mu}  U^{2-2\, \gamma}  +  \chi \gamma \left(\frac{4\chi (2-\gamma)}{\mu}\right)^{(2-\gamma)/\gamma}  }{ \frac{\nu}{C_{PF}^2}} \leq C (\text{data}).
 \end{gathered}
\end{equation}
\end{proof}
\begin{remark}
The stability analysis above assumes that $I_H$ is an
$L^2$-orthogonal projection and does not require the additional
conditions on $H$ and $\chi$ used in the analysis for a general
interpolation operator. For the corresponding assumptions and
stability analysis for a general interpolation operator, see
Appendix~\ref{app:general-interpolant}.
\end{remark}

\section{Error estimates}\label{error-estimates}
In this section, we derive an error estimate for the space-discrete, time-continuous scheme stated in (\ref{semifNLDANSE}).
We are interested in the error behavior up to a tolerance $\varepsilon$, where $0<\varepsilon\ll 1$.
Using the approximation properties of the finite element spaces stated in~\eqref{prop}, we assume that
\begin{equation}\label{h_assump_epsilon}
h^m \le C_3\,\varepsilon, \text{ and } h^{m+1}\leq \frac{|\bfu|_{m+1}}{2C_4} \epsilon,
\end{equation}
where $C_3>0$ is a constant independent of $h$ and $\varepsilon$.

 \begin{thm}
Let $I_H$ be the $L^2$-orthogonal projection defined in \eqref{l2-project}, and assume that \eqref{h_assump_epsilon} holds. Assume further that  $\frac{\nu}{4}-C_1^2 \mu H^2\geq 0$, $\frac{\nu}{4} - 2^{1-\gamma}C_1^2H\chi \geq 0$, $\chi \leq \mu$, and that
\begin{equation}\label{H_for_nonlinear}
    H\leq \left(\frac{C_1^2\alpha}{\nu}\right)^{\frac{\gamma}{2-2\gamma}} \mu^{\frac{\gamma}{2-2\gamma}}\epsilon^{\frac{2\gamma}{2-2\gamma}}.
\end{equation}
Then the error $\|\bfu-\bfv^h\|$ decays to the tolerance $\epsilon$ in finite time.
 \end{thm}
 
\begin{proof}
Define the error $e:= \bfu-\bfv^h$. We decompose the error into two pieces: $e= \eta -\phi^h$, where $\eta = \bfu- \Tilde{\bfu}$, and $\phi^h = \bfv^h-\Tilde{\bfu} \in \bfV^h$, $\Tilde{\bfu}\in \bfV^h$. With the assumption (\ref{h_assump_epsilon}), we have $\|\eta\|\leq \frac{\epsilon}{2}$. 
We take the inner product of (\ref{NSE}) and $
\Theta^h \in \bfV^h$ and subtract (\ref{semifNLDANSE}).  Setting $\Theta^h=\phi^h$, it yields
\begin{equation}
\begin{gathered}
\frac{1}{2}\frac{d \| \phi^h\|^2}{d t} + \nu \|\nabla \phi^h\|^2 +\mu (I_H(\bfv^h-\bfu), \phi^h) +\chi \left(\frac{I_H(\bfv^h-\bfu)}{\|I_H(\bfv^h-\bfu)\|^\gamma}, \phi^h\right) +(p,\nabla \cdot \phi^h) \\
\leq b(\bfu, \bfu,\phi^h) -b(\bfv^h, \bfv^h, \phi^h) +(\eta_t, \phi^h) + \nu (\eta, \phi^h).
\end{gathered}
\end{equation}
Consider the nudging term $\mu(I_H(\bfv^h-\bfu), \phi^h)$:
\begin{equation}
\begin{split}
\mu(I_H(\bfv^h-\bfu), \phi^h)
&=\mu \|I_H(\phi^h)\|^2 -\mu (I_H(\eta),I_H(\phi^h))\\
&\geq \mu \|I_H(\phi^h)\|^2- \frac{\mu}{2}\|I_H(\eta)\|^2 - \frac{\mu}{2}\|I_H(\phi^h)\|^2\\
&\geq \frac{\mu}{2}\|I_H(\phi^h)\|^2- \frac{\mu}{2}\|I_H(\eta)\|^2.
\end{split}
\end{equation}
We bound the nonlinear nudging term $\chi \left(\frac{I_H(\bfv^h-\bfu)}{\|I_H(\bfv^h-\bfu)\|^\gamma}, \phi^h\right)$ as the following:
\begin{equation}
\begin{gathered}
\chi \left(\frac{I_H(\bfv^h-\bfu)}{\|I_H(\bfv^h-\bfu)\|^\gamma}, \phi^h\right) = \chi \left(\frac{I_H(\phi^h-\eta)}{\|I_H(\phi^h-\eta)\|^\gamma}, \phi^h-\eta\right) + \chi \left(\frac{I_H(\phi^h-\eta)}{\|I_H(\phi^h-\eta)\|^\gamma}, \eta\right)\\
    = \chi \|I_H(\phi^h-\eta)\|^{2-\gamma} + \chi \left(\frac{I_H(\phi^h-\eta)}{\|I_H(\phi^h-\eta)\|^\gamma}, \eta\right)\\
    \geq \chi \|I_H(\phi^h-\eta)\|^{2-\gamma} -\chi \|I_H(\phi^h-\eta)\|^{1-\gamma}\|\eta\|\\
    \geq \chi \|I_H(\phi^h-\eta)\|^{2-\gamma} -\frac{1}{2} \chi \|I_H(\phi^h-\eta)\|^{2-\gamma} - \frac{\left( 2(1-\gamma)\right)^{1-\gamma}}{(2-\gamma)^{2-\gamma}} \chi \|\eta\|^{2-\gamma}\\
    =\frac{1}{2}\chi \|I_H(\phi^h-\eta)\|^{2-\gamma} -\frac{\left( 2(1-\gamma)\right)^{1-\gamma}}{(2-\gamma)^{2-\gamma}} \chi \|\eta\|^{2-\gamma} \\
    \geq \frac{1}{2}\chi 2^{1-\gamma}\|I_H(\phi^h)\|^{2-\gamma} - \frac{1}{2}\chi \|I_H(\eta)\|^{2-\gamma} -\frac{\left( 2(1-\gamma)\right)^{1-\gamma}}{(2-\gamma)^{2-\gamma}}\chi \|\eta\|^{2-\gamma}\\
    \geq  \frac{1}{2}\chi 2^{1-\gamma} \|I_H(\phi^h)\|^{2-\gamma}- \frac{1}{2}\chi C_I^{2-\gamma} \|\eta\|^{2-\gamma} -\frac{\left( 2(1-\gamma)\right)^{1-\gamma}}{(2-\gamma)^{2-\gamma}}\chi \|\eta\|^{2-\gamma}.
   \end{gathered} 
\end{equation}
The Navier--Stokes equations terms are treated as usual by
\begin{equation}
(\eta_t, \phi^h)\leq \|\eta_t\|_{-1}\|\nabla \phi^h\|\leq \frac{\nu}{12}\|\nabla \phi^h\|^2 +\frac{3}{\nu}\|\eta_t\|^2_{-1}.
\end{equation}
\begin{equation*}
\nu (\nabla \eta, \nabla \phi^h)\leq \nu\|\nabla \eta\|\|\nabla \phi^h\|\leq \frac{\nu}{12} \|\nabla \phi^h\|^2 + 3\nu \|\nabla \eta\|^2, \text{ and }
\end{equation*}
\begin{equation}
\begin{gathered}
(p,\nabla \cdot \phi^{h})=(p-q^h, \nabla \cdot \phi^h)\\
\leq \|p-q^h\|\|\nabla \cdot \phi^h\|\leq \frac{\nu}{12}\|\nabla \phi^h\|^2 +\frac{3}{\nu}\|p-q^h\|^2.
\end{gathered}
\end{equation}
The nonlinear terms $b(\bfu,\bfu,\phi^h)- b(\bfv^h,\bfv^h,\phi^h)$ satisfy
\begin{equation}
    \begin{gathered}
    b(\bfu, \bfu, \phi^h)- b(\bfv^h, \bfv^h, \phi^h)\\
    =b(\bfu,\bfu,\phi^h)+b(\phi^h-\bfv^h,\bfv^h,\phi^h)-b(\phi^h, \bfv^h,\phi^h)\\
        =b(\bfu,\bfu,\phi^h) + b(\eta -
        \bfu, \bfv^h, \phi^h) -b(\phi^h, \bfv^h,\phi^h)\\
        =b(\bfu, \eta,\phi^h) + b(\eta, \bfv^h,\phi^h) -b(\phi^h, \bfv^h, \phi^h).
    \end{gathered}
\end{equation}
Thus, term by term, 
\begin{equation}
\begin{split}
    b(\bfu,\eta,\phi^h)
    \leq K \|\nabla \eta\| \|\nabla \bfu\|\|\nabla \phi^h\| \leq \frac{\nu}{12}\| \nabla \phi^h\|^2+ \frac{3}{\nu} K^2\|\nabla \bfu\|^2 \|\nabla \eta\|^2.
\end{split}
\end{equation}
\begin{equation}
\begin{gathered}
    b(\eta, \bfv^h,\phi^h)\leq  K\|\nabla \eta\| \|\nabla \bfv^h\|\|\nabla \phi^h\|\\
    \leq \frac{\nu}{12} \|\nabla \phi^h\|^2 +\frac{3}{\nu}K^2 \|\nabla \bfv^h\|^2\|\nabla \eta\|^2.
\end{gathered}
\end{equation}
By Lemma \ref{nonlinear_bound_John}, we have
\begin{equation}
\begin{gathered}
    b(\phi^h, \bfv^h,\phi^h)\leq C_2\|\phi^h\|^{1/2}\|\nabla \phi^h\|^{3/2}\|\nabla \bfv^h\|\\
    \leq \frac{\nu}{12}\|\nabla \phi^h\|^2+C_2^4\frac{2187}{4}\nu^{-3}\|\nabla \bfv^h\|^4\|\phi^h\|^2.
\end{gathered}
\end{equation}
Combining the above, multiplying by $2$ and simplifying gives
\begin{equation}\label{error_put_together}
\begin{gathered}
\frac{d\|\phi^h\|^2}{dt} +\nu\|\nabla \phi^h\|^2  + \frac{\mu}{2}\|I_H(\phi^h)\|^2+ \chi 2^{1-\gamma}\|I_H(\phi^h)\|^{2-\gamma} \\
\leq \mu\|I_H(\eta)\|^2 + \frac{1}{2}\chi C_I^{2-\gamma} \|\eta\|^{2-\gamma}+ 2\frac{\left( 2(1-\gamma)\right)^{1-\gamma}}{(2-\gamma)^{2-\gamma}}\chi\|\eta\|^{2-\gamma}\\
+ \frac{6}{\nu}\|\eta_t\|^2_{-1} + \left(6\nu+ \frac{6}{\nu}K^2 \|\nabla \bfv^h\|^2+ \frac{6}{\nu}K^2\|\nabla \bfu\|^2 \right) \|\nabla \eta\|^2 \\
+C_2^4\frac{2187}{2}\nu^{-3} \|\nabla \bfv^h\|^4 \|\phi^h\|^2 + \frac{6}{\nu}\|p-q^h\|^2.
\end{gathered}
\end{equation}
Since $I_H$ is the $L^2$ orthogonal projection, we have
\begin{equation}
\|I_H(\phi^h)\|^2
=
\|\phi^h\|^2
-
\|(I-I_H)(\phi^h)\|^2.
\end{equation}
By the subadditivity of fractional powers, we have
\begin{equation}
    \|I_H(\phi^h)\|^{2-\gamma}\geq  \|\phi^h\|^{2-\gamma}-\|(I-I_H)(\phi^h)\|^{2-\gamma}.
\end{equation}
Then by (\ref{Ineq:Intrpolate}), we have
\begin{equation}
    \|I_H(\phi^h)\|^{2-\gamma}\geq  \|\phi^h\|^{2-\gamma}- C_1^{2-\gamma} H^{2-\gamma} \|\nabla \phi^h\|^{2-\gamma}.
\end{equation}
Using \eqref{Ineq:Intrpolate}, we obtain
\begin{equation}
\|I_H(\phi^h)\|^2
\ge
\|\phi^h\|^2
-
C_1^2 H^2 \|\nabla \phi^h\|^2.
\end{equation}
Thus,
\begin{equation}
\begin{gathered}
\frac{d \|\phi^h\|^2}{dt} + \frac{1}{2}\left(\nu-C_1^2 \mu H^2\right)\|\nabla \phi^h\|^2- 2^{1-\gamma}\chi  C_1^{2-\gamma} H^{2-\gamma} \|\nabla \phi^h\|^{2-\gamma}\\+\frac{\nu}{2}\|\nabla \phi^h\|^2 
+  \left(\frac{\mu}{2} - C_2^4\frac{2187}{2}\nu^{-3}\|\nabla \bfv^h\|^4 \right) \|\phi^h\|^2 
+\chi  2^{1-\gamma} \|\phi^h\|^{2-\gamma}
\\
 \leq \mu\|I_H(\eta)\|^2 
+ \frac{6}{\nu}\|\eta_t\|^2_{-1}+ \frac{1}{2}\chi C_I^{2-\gamma} \|\eta\|^{2-\gamma}+2\frac{\left( 2(1-\gamma)\right)^{1-\gamma}}{(2-\gamma)^{2-\gamma}}\chi \|\eta\|^{2-\gamma}\\ + \left(6\nu + \frac{6}{\nu}K^2\|\nabla \bfv^h\|^2+\frac{6}{\nu}K^2\|\nabla \bfu\|^2 \right) \|\nabla \eta\|^2 
+\frac{6}{\nu}\|p-q^h\|^2.\label{after-minkowski}
\end{gathered}
\end{equation}
Applying Young's inequality to $2^{1-\gamma}\chi C_1^{2-\gamma} H^{2-\gamma} \|\nabla \phi^h\|^{2-\gamma}$, we have
\begin{equation}
    2^{1-\gamma}\chi C_1^{2-\gamma} H^{2-\gamma} \|\nabla \phi^h\|^{2-\gamma}\leq \chi      2^{1-\gamma}\ C_1^2H\|\nabla \phi^h\|^2 +     2^{1-\gamma}\ \chi  H^{\frac{2-\gamma}{\gamma}}.
\end{equation}
For $\chi H^{\frac{2-\gamma}{\gamma}} $ term, since $\frac{1}{4}\nu -  2^{1-\gamma} C_1^2 H\chi \geq 0$,
\begin{equation}
        2^{1-\gamma}\chi H^{\frac{2-\gamma}{\gamma}} \leq \frac{\nu}{4 C_1^2} H^{\frac{2-2\gamma}{\gamma}}. 
\end{equation}
Since (\ref{H_for_nonlinear}) holds, we have 
\begin{equation}
\frac{\nu}{4 C_1^2} H^{\frac{2-2\gamma}{\gamma}}\leq \frac{\alpha}{4}\mu \|\phi^h\|^2. 
\end{equation}
Since $\frac{1}{4}\nu - 2^{1-\gamma} C_1^2 H\chi \geq 0$ and $\frac{\nu}{4}-C_1^2 \mu H^2\geq 0$,
(\ref{after-minkowski}) becomes
\begin{equation}\label{before_cases}
\begin{gathered}
\frac{d \|\phi^h\|^2}{dt}+\frac{\nu}{4}\|\nabla \phi^h\|^2 
+ \frac{\alpha}{4} \mu  \|\phi^h\|^2 
+2^{1-\gamma}\chi  \|\phi^h\|^{2-\gamma}
\\
 \leq C_I^2\mu\|\eta\|^2 
+ \frac{6}{\nu}\|\eta_t\|^2_{-1} + \left(6\nu + \frac{6}{\nu}K^2\|\nabla \bfv^h\|^2+\frac{6}{\nu}K^2\|\nabla \bfu\|^2 \right) \|\nabla \eta\|^2\\
+\frac{6}{\nu}\|p-q^h\|^2 + \chi \left(\frac{1}{2} C_I^{2-\gamma} +2\frac{\left( 2(1-\gamma)\right)^{1-\gamma}}{(2-\gamma)^{2-\gamma}}\right)  \|\eta\|^{2-\gamma}. 
\end{gathered}
\end{equation}
We will consider two cases. \textit{When $\|\phi^h\|\geq 1$, $2^{1-\gamma}\chi \|\phi^h\|^{2-\gamma}$ can be absorbed and is therefore omitted in (\ref{before_cases}).} Multiply both sides by the integration factor $\exp\big[\frac{\alpha}{4}\mu t\big]$. We integrate over time from $0$ to $t$, and multiply $\exp \big[-\frac{\alpha}{4}\mu t\big]$ to both sides, (\ref{before_cases}) becomes
\begin{equation*}
\begin{gathered}
\|\phi^h(t)\|^2 + \frac{\nu}{4}\int_{0}^t \exp\bigg[-\frac{\alpha}{4}\mu (t-t')\bigg]\|\nabla \phi^h\|^2\, dt'  \leq \exp\bigg[-\frac{\alpha}{4} \mu t\bigg] \|\phi^h(0)\|^2\\
+ \int_{0}^t \exp\bigg[-\frac{\alpha}{4} \mu (t-t')\bigg]
C_I^2\mu\|\eta\|^2\, d t' \\
+ \frac{6}{\nu}\|\eta_t\|^2_{-1}
+ \left(6\nu + \frac{6}{\nu}K^2\|\nabla \bfv^h\|^2+\frac{6}{\nu}K^2\|\nabla \bfu\|^2 \right) \|\nabla \eta\|^2
+\frac{6}{\nu}\|p-q^h\|^2 \\
+\int_{0}^t  \exp\bigg[-\frac{\alpha}{4} \mu (t-t')\bigg] \chi \left(\frac{C_I^{2-\gamma}}{2} +2\frac{\left( 2(1-\gamma)\right)^{1-\gamma}}{(2-\gamma)^{2-\gamma}}\right) \|\eta\|^{2-\gamma}\, dt'. 
\end{gathered}
\end{equation*}
We bound 
\begin{equation}
    \int_{0}^t \exp\bigg[-\frac{\alpha}{4} \mu (t-t')\bigg]
C_I^2\mu\|\eta\|^2\, d t'
\end{equation}
as follows.
\begin{equation*}
\begin{gathered}
\int_{0}^T \exp\bigg[-\frac{\alpha}{4} \mu (t-t')\bigg]
C_I^2\mu\|\eta\|^2\, d t'
\leq \max_{0\leq t \leq T} \|\eta\|^2 C_I^2 \mu \exp\bigg[-\frac{\alpha}{4} \mu t\bigg]\int_{0}^t \exp\bigg[\frac{\alpha}{4} \mu t'\bigg] \, dt'\\
=\max_{0\leq t \leq T} \|\eta\|^2 C_I^2 \mu \exp\bigg[-\frac{\alpha}{4} \mu t\bigg] \left(\exp\bigg[\frac{\alpha}{4} \mu t\bigg] -1\right)\frac{4}{\alpha} \mu^{-1}\\
=\max_{0\leq t \leq T} \|\eta\|^2 C_I^2 \frac{4}{\alpha} \left(1-\exp\bigg[-\frac{\alpha}{4} \mu t\bigg]\right)
\leq  C_I^2 \frac{4}{\alpha}\max_{0\leq t \leq T} \|\eta\|^2. 
\end{gathered}
\end{equation*}
Similarly, under the assumption $\chi \leq \mu$, we have
\begin{equation}
\begin{gathered}
\int_{0}^t  \exp\bigg[-\frac{\alpha}{4} \mu (t-t')\bigg] \chi \left(\frac{1}{2} C_I^{2-\gamma} +2\frac{\left( 2(1-\gamma)\right)^{1-\gamma}}{(2-\gamma)^{2-\gamma}}\right) \|\eta\|^{2-\gamma}\, dt'\\
\leq \left( \frac{1}{2} C_I^{2-\gamma} +2\frac{\left( 2(1-\gamma)\right)^{1-\gamma}}{(2-\gamma)^{2-\gamma}}\right) \frac{4}{\alpha} \frac{\chi}{\mu} \max_{0\leq t\leq T} \|\eta\|^{2-\gamma}\\
\leq \left( \frac{1}{2} C_I^{2-\gamma} +2\frac{\left( 2(1-\gamma)\right)^{1-\gamma}}{(2-\gamma)^{2-\gamma}}\right)  \frac{4}{\alpha} \max_{0\leq t\leq T} \|\eta\|^{2-\gamma}.
\end{gathered}
\end{equation}
We simplify by the following estimates.
\begin{equation}\label{estimate_1}
    \begin{gathered}
        \int_{0}^T \|\nabla \bfu\|^2 \|\nabla \eta\|^2\, dt \leq \|\nabla \bfu\|^2_{L^4(0,T;L^2)}\|\nabla \eta\|^2_{L^4(0,T;L^2)},
    \end{gathered}
\end{equation}

\begin{equation}\label{estimate_2}
    \begin{gathered}
        \int_{0}^T \|\nabla \bfv^h\|^2 \|\nabla \eta\|^2\, dt \leq \|\nabla \bfv^h\|^2_{L^4(0,T;L^2)}\|\nabla \eta\|^2_{L^4(0,T;L^2)}.
    \end{gathered}
\end{equation}
Hence, we have
\begin{equation}\label{before_triangle}
\begin{gathered}
\|\phi^h(t)\|^2 + \frac{\nu}{4}\int_{0}^t \exp\bigg[-\frac{\alpha}{4}\mu (t-t')\bigg]\|\nabla \phi^h\|^2\, dt'  \leq \exp\bigg[-\frac{\alpha}{4} \mu t\bigg] \|\phi^h(0)\|^2\\
+ C_I^2 \frac{4}{\alpha}\max_{0\leq t \leq T} \|\eta\|^2 + \left( \frac{1}{2} C_I^{2-\gamma} +2\frac{\left( 2(1-\gamma)\right)^{1-\gamma}}{(2-\gamma)^{2-\gamma}}\right) \frac{4}{\alpha} \max_{0\leq t\leq T} \|\eta\|^{2-\gamma}
+ \frac{6}{\nu} \int_{0}^t \|\eta_t\|^2_{-1}\, dt'\\
+  \frac{6}{\nu}K^2\left(\|\nabla \bfv^h\|^2_{L^4(0,T;L^2)}+\|\nabla \bfu\|^2_{L^4(0,T;L^2)} \right) \|\nabla \eta\|^2_{L^4(0,T;L^2)}\\
+6\nu \|\nabla \eta\|^2_{L^2(0,T; L^2)}
+\frac{6}{\nu}\int_{0}^t\|p-q^h\|^2\, dt'. 
\end{gathered}
\end{equation}
Using the triangle inequality $\|e\| \leq \|\eta\|+\|\phi^h\|$, it gives
\begin{equation}\label{triangle-eqn}
\begin{gathered}
\|e(t)\|^2 + \frac{\nu}{4}\int_{0}^t \exp\bigg[-\frac{\alpha}{4}\mu (t-t')\bigg]\|\nabla e\|^2\, dt'\leq \|\eta(t)\|^2 \\
+ \nu \int_{0}^t \|\nabla \eta\|^2\, dt' + \|\phi^h(t)\|^2 + \frac{\nu}{4}\int_{0}^t \exp\bigg[-\frac{\alpha}{4}\mu (t-t')\bigg]\|\nabla \phi^h\|^2\, dt'. 
\end{gathered}
\end{equation}
Thus, combining (\ref{triangle-eqn}) with (\ref{before_triangle}), we have
\begin{equation}
\begin{gathered}
 \|e(t)\|^2 + \frac{\nu}{4}\int_{0}^t \exp\bigg[-\frac{\alpha}{4}\mu (t-t')\bigg]\|\nabla e\|^2\, dt'\leq 
 \exp\bigg[-\frac{\alpha}{4} \mu t\bigg] \|e(0)\|^2\\
+ \left(C_I^2 \frac{4}{\alpha}+2 \right)\max_{0\leq t \leq T} \|\eta\|^2 
+ \left( \frac{1}{2} C_I^{2-\gamma} +2\frac{\left( 2(1-\gamma)\right)^{1-\gamma}}{(2-\gamma)^{2-\gamma}}\right) \frac{4}{\alpha} \max_{0\leq t\leq T} \|\eta\|^{2-\gamma}
\\+ \frac{6}{\nu} \int_{0}^t \|\eta_t\|^2_{-1}\, dt'
+  \frac{6}{\nu}K^2\left(\|\nabla \bfv^h\|^2_{L^4(0,T;L^2)}+\|\nabla \bfu\|^2_{L^4(0,T;L^2)} \right) \|\nabla \eta\|^2_{L^4(0,T;L^2)}\\
+7\nu \|\nabla \eta\|^2_{L^2(0,T; L^2)}
+\frac{6}{\nu}\int_{0}^t\|p-q^h\|^2\, dt'. 
\end{gathered}
\end{equation}
Thus, by the definition of $\eta$, we have
\begin{equation}
\begin{gathered}
 \|e(t)\|^2 + \frac{\nu}{4}\int_{0}^t \exp\bigg[-\frac{\alpha}{4}\mu (t-t')\bigg]\|\nabla e\|^2\, dt'\leq   
 \exp\bigg[-\frac{\alpha}{4} \mu t\bigg] \|e(0)\|^2\\
+ \left(C_I^2 \frac{4}{\alpha}+2 \right)\max_{0\leq t \leq T} \|\bfu-\tilde{\bfu}\|^2 + \left( \frac{1}{2} C_I^{2-\gamma} +2\frac{\left( 2(1-\gamma)\right)^{1-\gamma}}{(2-\gamma)^{2-\gamma}}\right) \frac{4}{\alpha} \max_{0\leq t\leq T} \|\bfu-\tilde{\bfu}\|^{2-\gamma}\\
+ \frac{6}{\nu} \int_{0}^t \|(\bfu-\tilde{\bfu})_t\|^2_{-1}\, dt'
+  \frac{6}{\nu}K^2\left(\|\nabla \bfv^h\|^2_{L^4(0,T;L^2)}+\|\nabla \bfu\|^2_{L^4(0,T;L^2)} \right) \|\nabla (\bfu-\tilde{\bfu})|^2_{L^4(0,T;L^2)}\\
+7\nu \|\nabla (\bfu-\tilde{\bfu})\|^2_{L^2(0,T; L^2)}
+\frac{6}{\nu}\int_{0}^t\|p-q^h\|^2\, dt'. 
\end{gathered}
\end{equation}
\textit{When $\|\phi^h\|<1$, the nonlinear term $2^{1-\gamma}\chi \|\phi^h\|^{2-\gamma}$ dominates the quadratic terms. Therefore, we drop $\frac{\alpha}{4} \mu \|\phi^h\|^2$ and $\frac{\nu}{4}\|\nabla \phi^h\|^2$, thus (\ref{before_cases}) becomes}
\begin{equation}
\label{non-linear}
\begin{gathered}
\frac{d \|\phi^h\|^2}{dt} +2^{1-\gamma}\chi  \|\phi^h\|^{2-\gamma}
\leq C_I^2\mu\|\eta\|^2 + \frac{6}{\nu}\|\eta_t\|^2_{-1} \\
+ \left(6\nu + \frac{6}{\nu}K^2\|\nabla \bfv^h\|^2+\frac{6}{\nu} K^2\|\nabla \bfu\|^2 \right) \|\nabla \eta\|^2
+\frac{6}{\nu}\|p-q^h\|^2\\
+ \chi \left(\frac{1}{2} C_I^{2-\gamma} +2\frac{\left( 2(1-\gamma)\right)^{1-\gamma}}{(2-\gamma)^{2-\gamma}}\right)\|\eta\|^{2-\gamma}. 
\end{gathered}
\end{equation}
We define 
\begin{equation}
\begin{gathered}
G:= C_I^2\mu\|\eta\|^2 + \frac{6}{\nu}\|\eta_t\|^2_{-1} + \left(6\nu + \frac{6}{\nu}K^2\|\nabla \bfv^h\|^2+\frac{6}{\nu} K^2\|\nabla \bfu\|^2 \right) \|\nabla \eta\|^2\\
+\frac{6}{\nu}\|p-q^h\|^2 + \chi \left( \frac{1}{2} C_I^{2-\gamma} +2\frac{\left( 2(1-\gamma)\right)^{1-\gamma}}{(2-\gamma)^{2-\gamma}}\right)\|\eta\|^{2-\gamma}.
\end{gathered}
\end{equation}
Thus, we can rewrite $(\ref{non-linear})$ as 
\begin{equation}
    \frac{d \|\phi^h\|^2}{dt} +2^{1-\gamma}\chi  \|\phi^h\|^{2-\gamma} \leq G.
\end{equation}
Divide both sides by $\|\phi^h\|^{2-\gamma}$, we have
\begin{equation}
    \frac{d \|\phi^h\|^2}{dt} \|\phi^h\|^{\gamma -2} + 2^{1-\gamma} \chi \leq G \frac{1}{\|\phi^h\|^{2-\gamma}}.
\end{equation} 
Since
\begin{equation}
    \|\phi^h\|= \|\eta -e \|\geq \|e\| - \|\eta\| \geq \frac{\epsilon}{2},
\end{equation}
we have
\begin{equation}
    \frac{d \|\phi^h\|^2}{dt} (\|\phi^h\|^2)^{\gamma/2 -1} + 2^{1-\gamma} \chi \leq G \left( \frac{2}{\epsilon}\right)^{2-\gamma}.
\end{equation}
Integrating over time from $0$ to $t$,
\begin{equation}
\frac{2}{\gamma}(\|\phi^h(t)\|^2)^{\gamma/2} - \frac{2}{\gamma}(\|\phi^h(0)\|^2)^{\gamma/2} + 2^{1-\gamma}\chi t \leq  \left( \frac{2}{\epsilon}\right)^{2-\gamma} \int_{0}^t G(t')\, dt'.
\end{equation}
Dividing both sides by $\frac{2}{\gamma}$ and rearranging terms, we have
\begin{equation}(\|\phi^h(t)\|^2)^{\gamma/2}\leq (\|\phi^h(0)\|^2)^{\gamma/2} - \frac{\gamma}{2^{\gamma}}  \chi t + \frac{\gamma}{2}\left( \frac{2}{\epsilon}\right)^{2-\gamma} \int_{0}^t G(t')\, dt'.
\end{equation}
Now, we estimate $\frac{\gamma}{2}\left(\frac{2}{\epsilon}\right)^{2-\gamma} \int_{0}^t G(t')\, dt'$.
\begin{equation}
\begin{gathered}
    \frac{\gamma}{2}\left(\frac{2}{\epsilon}\right)^{2-\gamma} \int_{0}^t G(t')\, dt' \leq \frac{\gamma}{2}\left(\frac{2}{\epsilon}\right)^{2-\gamma}\Bigg[\left(C_I^2 \frac{4}{\alpha}+2 \right)\max_{0\leq t \leq T} \|\eta\|^2 \\
    + \left( \frac{1}{2} C_I^{2-\gamma} +2\frac{\left( 2(1-\gamma)\right)^{1-\gamma}}{(2-\gamma)^{2-\gamma}}\right) \frac{4}{\alpha} \max_{0\leq t\leq T} \|\eta\|^{2-\gamma}
+ \frac{6}{\nu} \int_{0}^t \|\eta_t\|^2_{-1}\, dt'\\
+  \frac{6}{\nu}K^2\left(\|\nabla \bfv^h\|^2_{L^4(0,T;L^2)}+\|\nabla \bfu\|^2_{L^4(0,T;L^2)} \right) \|\nabla \eta\|^2_{L^4(0,T;L^2)}\\
+6\nu \|\nabla \eta\|^2_{L^2(0,T; L^2)}
+\frac{6}{\nu}\int_{0}^t\|p-q^h\|^2\, dt'\Bigg].
\end{gathered}
\end{equation}
With the assumption (\ref{h_assump_epsilon}),
\begin{equation}
\begin{gathered}
    \frac{\gamma}{2}\left(\frac{2}{\epsilon}\right)^{2-\gamma} \int_{0}^t G(t')\, dt' \leq C(\nu, C_I, C_4 \alpha, \gamma, K, \bfu, \bfv^h)\frac{\gamma}{2}2^{2-\gamma} \left(\frac{h^m}{\epsilon}\right)^{2-\gamma}\\
    \leq C(\nu, C_I, C_4, \alpha, \gamma, K, \bfu, \bfv^h).
\end{gathered}
\end{equation}
By the uniform bound obtained above, there exists a constant
$C_0\ge0$, independent of time, such that
\begin{equation}
    \frac{\gamma}{2}\left(\frac{2}{\epsilon}\right)^{2-\gamma}
    \int_0^t G(s)\,ds \le C_0.
\end{equation}
We define
\begin{equation}
    t^\star :=
    \frac{2^\gamma}{\gamma\chi}
    \left(
        \|\phi^h(0)\|^\gamma+C_0
        -\left(\frac{\epsilon}{2}\right)^\gamma
    \right).
\end{equation}
If the tolerance has not already been reached, the preceding
estimate gives $\|\phi^h(t^\star)\|\le\epsilon/2$.
By the triangle inequality,
\begin{equation}
    \|e(t^\star)\|
    \le \|\eta(t^\star)\|+\|\phi^h(t^\star)\|
    \le \frac{\epsilon}{2}+\frac{\epsilon}{2}
    =\epsilon.
\end{equation}
Thus, the error reaches the prescribed tolerance in finite time,
no later than $t^\star$.
\end{proof}

\section{Numerical Experiment}\label{numerical_test}

In this section, we test the nonlinear nudging method in four numerical experiments. In Section~\ref{S5:1}, we verify the temporal convergence rate using an analytical solution. In Section~\ref{S5;2}, we consider a forced two-dimensional flow in a periodic box, in which the nonlinear nudging algorithm was originally tested with pseudo-spectral methods \cite{CLT23}. In Sections~\ref{sec:flatplate} and~\ref{sec:cylinder}, we consider flows past an obstacle in a channel, with a thin flat plate and  a smooth circular obstacle, respectively, which involve physical boundary conditions and cannot be treated with spectral methods. In all experiments, we compare nonlinear nudging ($\chi>0$) with linear nudging ($\chi=0$), initialize the assimilated solution from rest, and evaluate the denominator of the nonlinear nudging term at the previous time level, as in the remark following Algorithm~\ref{Algorithm}. We also report the effective nudging coefficient $\mu+\chi\|I_H(\bfu_n-\bfv^h_n)\|^{-\gamma}$, which measures the total strength of the nudging at each time step.

\subsection{Accuracy and convergence tests using an analytical solution}\label{S5:1}
We validate the accuracy of the nonlinear data assimilation using a test with exact solutions. The spatial discretization employs a Scott–Vogelius finite element pair on a barycenter-refined mesh, and time integration is performed using the BDF2. We verify the expected second-order temporal convergence.

In the numerical test, we linearize the nonlinear data assimilation term:
\begin{equation}
    \chi \|I_H(\bfu_n)-I_H(\bfv^h_n)\|^{-\gamma} (I_H (\bfu_{n+1})-I_H(\bfv^h_{n+1})).
\end{equation}
The setup of the test is in a unit square domain $\Omega =(0,1)\times (0,1)$. The exact velocity and pressure are 
\begin{equation}
    \bfu(x,y,t) = \langle e^{t} \cos(y), e^{t} \sin(x) \rangle^\top, \quad \text{and } p(x,y,t) = (x-y) (1+t).
\end{equation}
We set the characteristic velocity $U=1$, the characteristic length scale $L=1$, and the kinematic viscosity $\nu=1$, so that the Reynolds number
\[
\mathcal{R}e = \frac{UL}{\nu} = 1.
\]
We then calculate the corresponding body force $f(x,y,t)$.  We generate a fine mesh with $55554$ degrees of freedom (DOF) so the spacial errors are small and negligible. We run the simulations in the time interval $[0,2]$. We set $\mu=10$, $\chi=10$ and $\gamma=0.5$. Timestep sizes $1, 1/2, 1/4, 1/8, 1/16$, and $1/32$. Table \ref{tab:error-main} demonstrates the convergence rate for nonlinear DA with $\gamma=0.5$. We observe a second order convergence rate.
\begin{table}[ht]
    \centering
        \begin{tabular}{|c|c|c|c|}
            \hline
            $\Delta t$ & $\|u-v\|$ & rate \\
            \hline
            1       & 0.0043   & -        \\
            1/2     & 0.0010   & 2.04    \\
            1/4     & 0.00030  & 1.76   \\
            1/8     & 8.3e-5   & 1.88    \\
            1/16    & 2.2e-5   & 1.94    \\
            1/32    & 5.6e-6   & 1.94   \\
            \hline
        \end{tabular}
    \caption{Errors, rates of convergence for the nonlinear DA. We obtain the second-order error rates in time.}
    \label{tab:error-main}
\end{table}

\subsection{Forced two-dimensional flow  in a periodic box}\label{S5;2}
In this section, we revisit the two-dimensional turbulence test of \cite{CLT23} Section 5, which was carried out with a pseudo-spectral method, and repeat it with the finite element method. The aim is to check whether the super-exponential convergence of nonlinear nudging persists under a finite element discretization and a finite element interpolant $I_H$.

\subsection*{Simulation Setup}
The domain is $\Omega=[0,2\pi)^2$ with periodic boundary conditions. The fine DNS  mesh is a uniform $128\times128$ grid of squares, each split into two triangles ($h=2\pi/128$), and we use Taylor--Hood $(P_2,P_1)$ elements, about $1.5\times 10^{5}$ degrees of freedom in total. For time integration, we use the second-order backward differentiation formula (BDF2) with $\Delta t = 5\times10^{-4}$. The convective term, in the skew-symmetric form $b(\cdot,\cdot,\cdot)$, is treated explicitly by extrapolation, $2b(\bfv^h_n,\bfv^h_n,\Theta^h)-b(\bfv^h_{n-1},\bfv^h_{n-1},\Theta^h)$, while the Stokes part and both nudging terms are implicit. As in the remark following Algorithm \ref{Algorithm}, the denominator of the nonlinear nudging term is evaluated at the previous time level, so that each time step requires only one linear solve. 

\begin{figure}
    \centering
    \begin{subfigure}{0.48\textwidth}
        \centering
        \includegraphics[width=\linewidth]{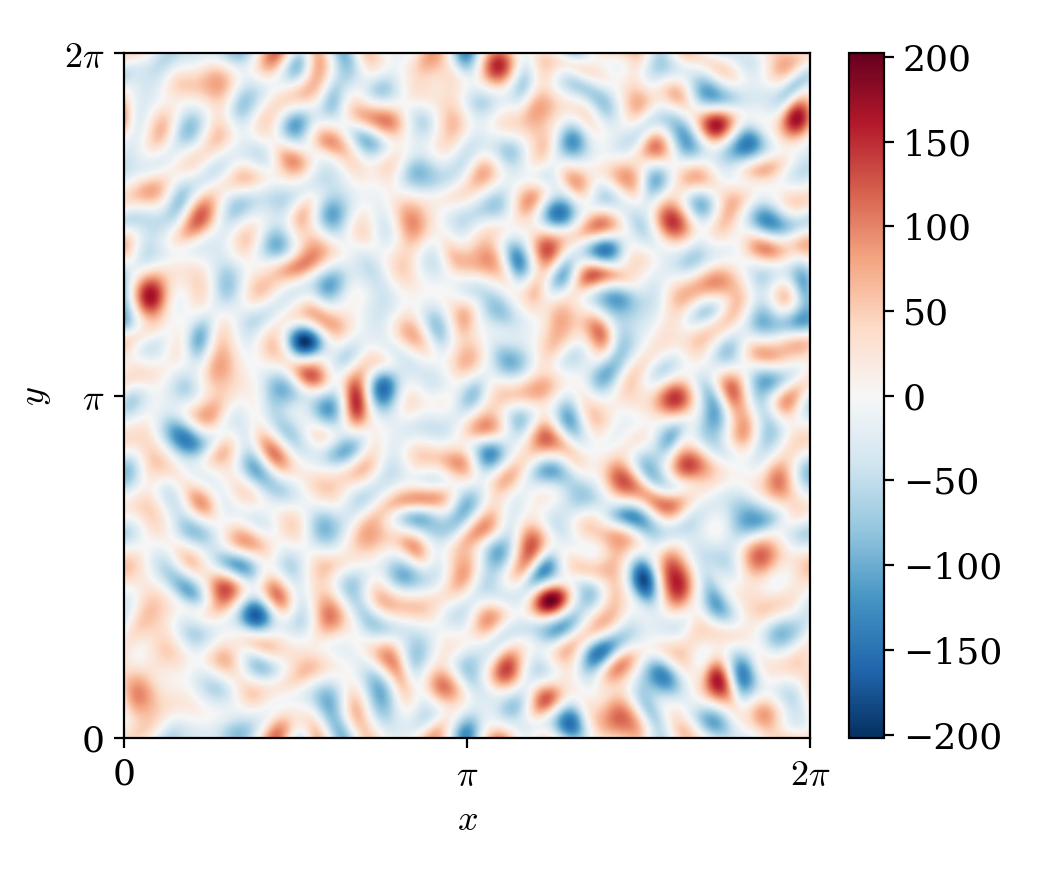}
        \caption{Vorticity of the initial condition $\bfu_0^h$.}
        \label{fig:N128_vort}
    \end{subfigure}
    \hfill
    \begin{subfigure}{0.48\textwidth}
        \centering
        \includegraphics[width=\linewidth]{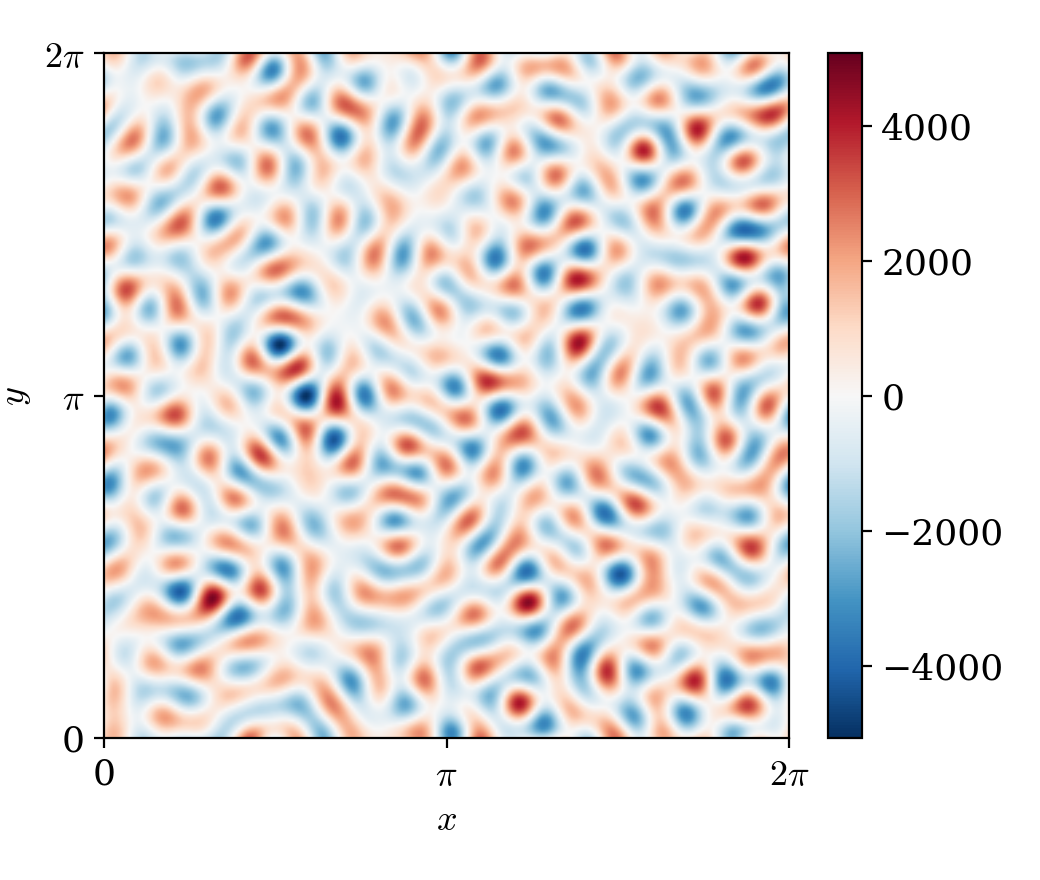}
        \caption{Curl of the body force, $\nabla\times\bff$.}
        \label{fig:N128_force}
    \end{subfigure}
    \begin{subfigure}{0.48\textwidth}
        \centering
        \includegraphics[width=\linewidth]{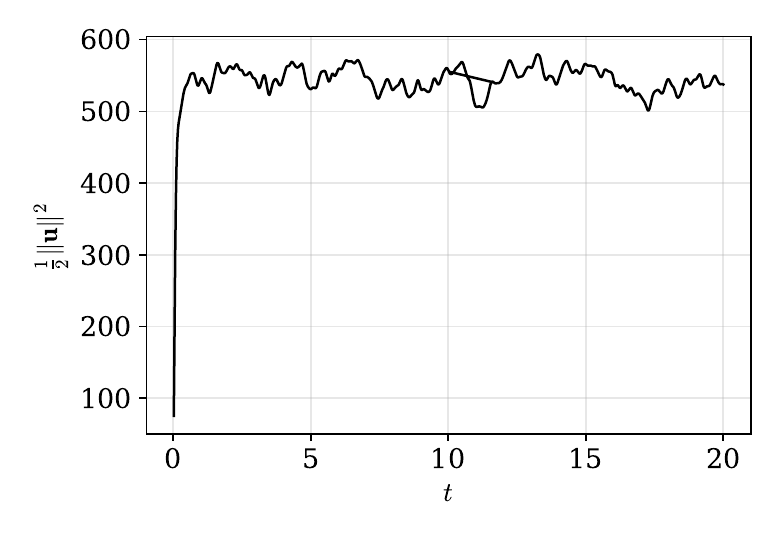}
        \caption{Energy $\frac12\|\bfu(t)\|^2$ during the spin-up.}
        \label{fig:N128_energy}
    \end{subfigure}
    \hfill
    \begin{subfigure}{0.48\textwidth}
        \centering
        \includegraphics[width=\linewidth]{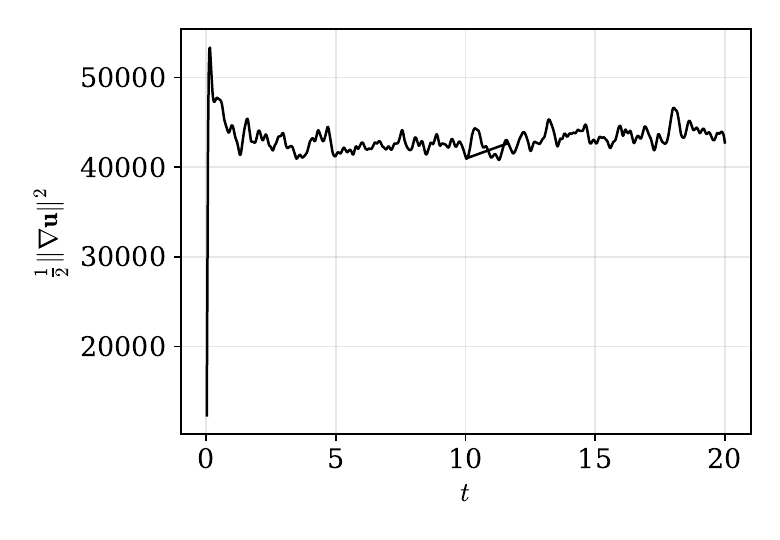}
        \caption{Enstrophy $\frac12\|\nabla\bfu(t)\|^2$ during the spin-up.}
        \label{fig:N128_enstrophy}
    \end{subfigure}
    \caption{Reference solution for the periodic box test, $128\times128$ mesh.}
    \label{fig:N128_truth}
\end{figure}

We set the kinematic viscosity to $\nu=0.128$. The body force is $\bff=\nabla^{\perp} g$, where $g$ has random Fourier coefficients on the shell $4\leq|\bk|\leq16$, normalized so that $\|\bff\|=640$. With this forcing the flow is time-dependent and statistically steady; see Figures \ref{fig:N128_energy} and \ref{fig:N128_enstrophy}. 

To generate the initial condition $\bfu_0^h$ for the reference solution, we simulate the system from rest up to $t=20$. The resulting flow is statistically steady, with energy $\frac12\|\bfu\|^2\approx 540$ and enstrophy $\frac12\|\nabla\bfu\|^2\approx 4.3\times10^4$; see Figure \ref{fig:N128_truth}.

The observations are given by local averages over a $16\times16$ grid of coarse squares, i.e., $H=2\pi/16=8h$. This interpolant is the $L^2$ orthogonal projection onto piecewise constants on the coarse mesh; it satisfies \eqref{Ineq:Intrpolate} and is the setting used in the error analysis of Section 4. The assimilated solution is initialized with $\bfv^h_0=\mathbf{0}$, so that the initial error equals $\|\bfu_0^h\|\approx 33$. We compare the following cases:
\begin{itemize}
    \item no data assimilation ($\mu=\chi=0$), 
    \item linear nudging with $\mu=2$ and with $\mu=4$ ($\chi=0$), 
    \item nonlinear nudging with $\mu=\chi=2$ and $\gamma=0.1,\ 0.25,\ 0.5$.
\end{itemize}
The linear case $\mu=4$ is included because it equals $\mu+\chi$, i.e., the effective coefficient of the nonlinear scheme when $\|I_H(\bfu_n-\bfv_n^h)\|=1$. This allows us to distinguish the contribution of the nonlinear feedback from that of a larger linear nudging parameter.

\subsection*{Results}
\begin{figure}
    \centering
    \begin{subfigure}{0.48\textwidth}
        \centering
        \includegraphics[width=\linewidth]{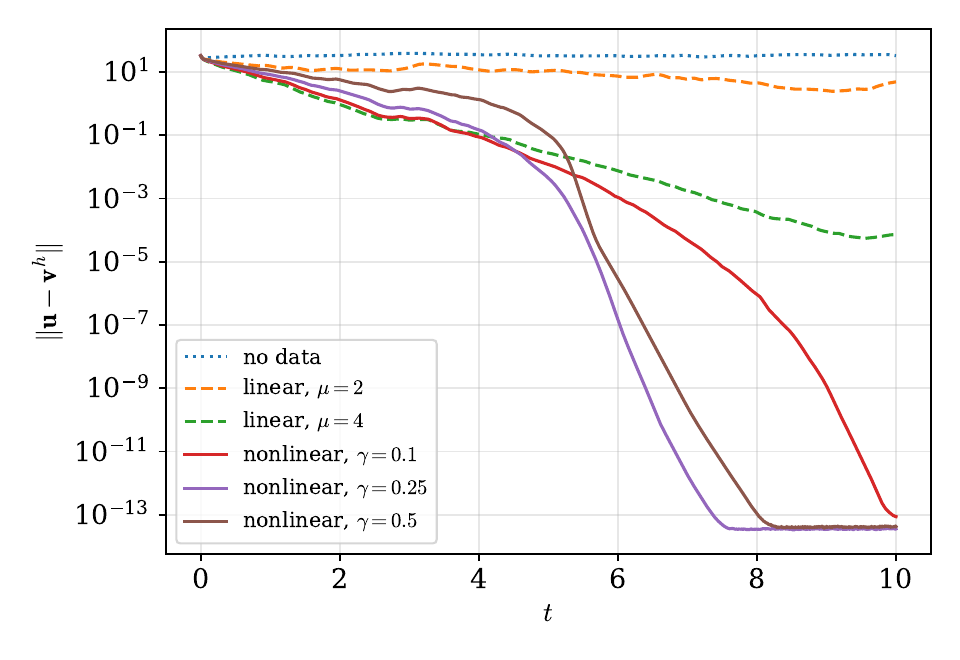}
        \caption{Error $\|\bfu-\bfv^h\|$.}
        \label{fig:N128_err}
    \end{subfigure}
    \hfill
    \begin{subfigure}{0.48\textwidth}
        \centering
        \includegraphics[width=\linewidth]{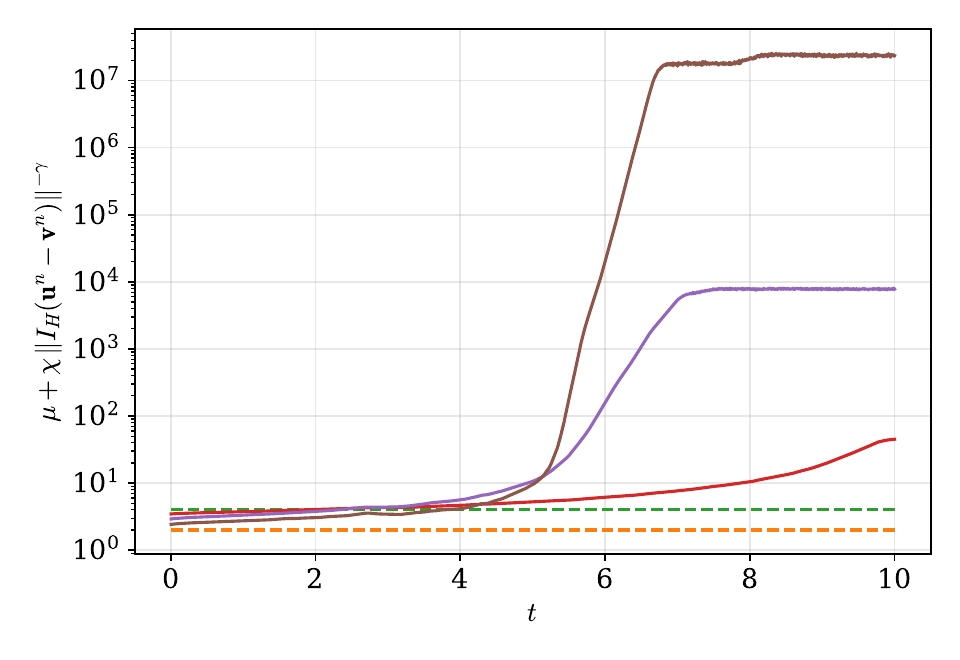}
        \caption{Effective coefficient $\mu+\chi\|I_H(\bfu_n-\bfv^h_n)\|^{-\gamma}$.}
        \label{fig:N128_mueff}
    \end{subfigure}
    \caption{Periodic box test: linear nudging ($\mu=2,4$, $\chi=0$) and nonlinear nudging ($\mu=\chi=2$, $\gamma=0.1,0.25,0.5$). All nonlinear cases reach machine precision, while the linear cases do not by $t=10$.}
    \label{fig:N128_DA}
\end{figure}

Figure \ref{fig:N128_DA} shows the error $\|\bfu_n-\bfv^h_n\|$ (Figure \ref{fig:N128_err}) and the effective nudging coefficient $\mu+\chi\|I_H(\bfu_n-\bfv^h_n)\|^{-\gamma}$ (Figure \ref{fig:N128_mueff}) for all cases. Without data assimilation, the error remains at about $30$, the size of the solution itself. Linear nudging with $\mu=2$ reduces the error only to about $5$ by $t=10$, while with $\mu=4$ the error decays exponentially and reaches $7.3\times10^{-5}$ at $t=10$.

All nonlinear cases reach machine precision, $\mathcal{O}(10^{-14})$, by $t=10$. Up to $t\approx4$, their effective coefficients remain between about $2$ and $6$, close to the linear value $\mu=4$, and their errors decay at a rate comparable to that of the linear case $\mu=4$. After $t\approx4$, the effective coefficients grow rapidly, and the errors decay faster than exponentially. For $\gamma=0.25$, the error drops from $10^{-6}$ at $t=6$ to $10^{-12}$ at $t=7$; for $\gamma=0.5$, from $10^{-6}$ at $t=6$ to $10^{-12}$ at $t=8$. Once the error reaches machine precision, the effective coefficient stops growing and levels off, at about $8\times10^{3}$ for $\gamma=0.25$ and $3\times10^{7}$ for $\gamma=0.5$. For $\gamma=0.1$, the effective coefficient grows more slowly, to about $45$ at $t=10$, and the error reaches $10^{-6}$ at $t=8.0$ and $10^{-12}$ at $t=9.7$.

The dependence on $\gamma$ is not monotone. At early times, a larger $\gamma$ gives a smaller effective coefficient and a larger error; for $t<5$, the case $\gamma=0.5$ has the largest error among the nonlinear cases. Afterwards, a larger $\gamma$ gives a faster growth of the effective coefficient and a steeper decay of the error. The cases $\gamma=0.25$ and $\gamma=0.5$ perform similarly, and both reach machine precision about $2$ time units before $\gamma=0.1$.

\subsection{Flow over a flat plate obstacle}\label{sec:flatplate}

\begin{figure}
    \centering
\includegraphics[width=0.5\linewidth]{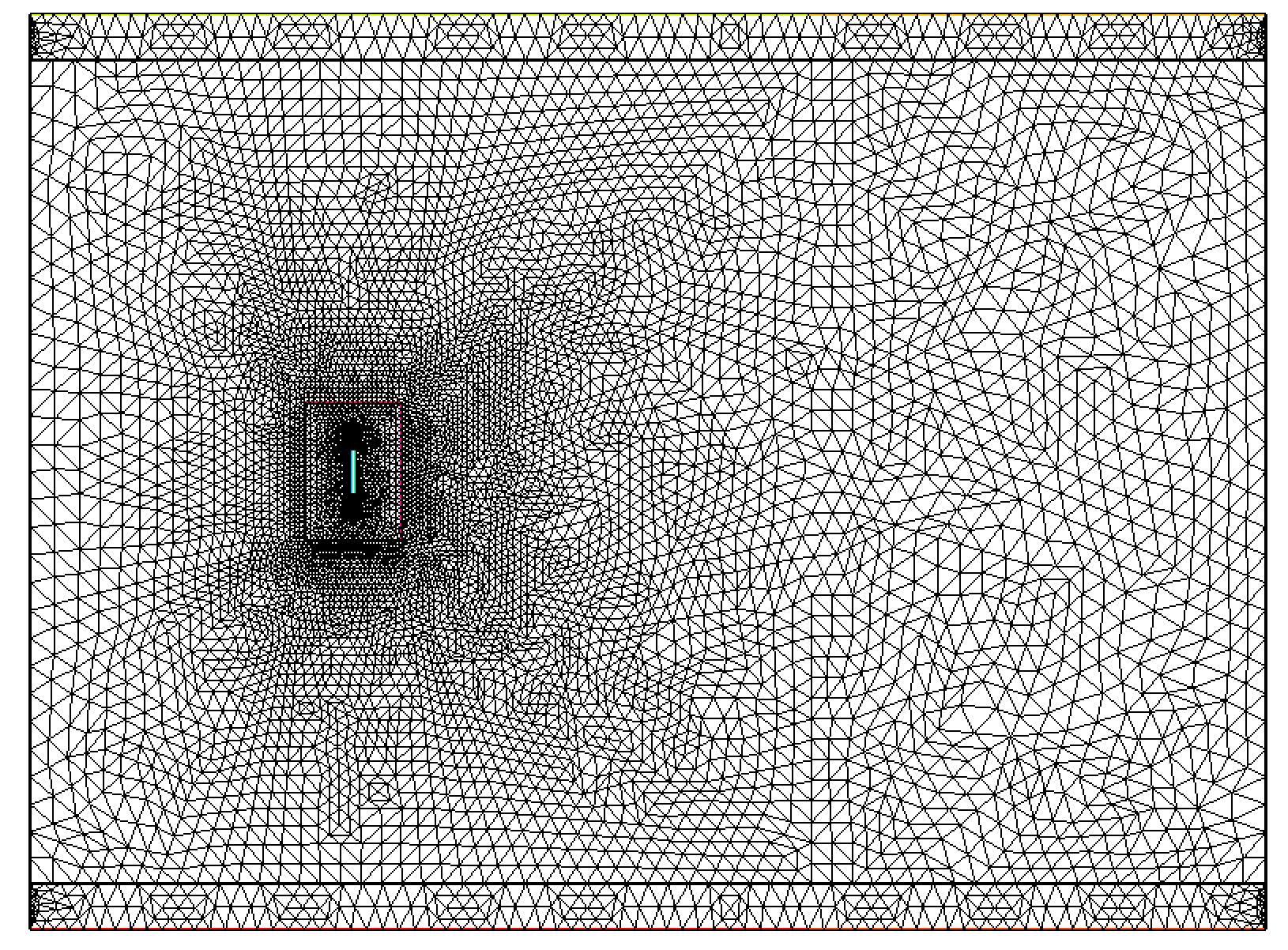}
    \caption{Domain with meshing.}
    \label{fig:test3-mesh}
\end{figure}

In this section, we test the nonlinear nudging method on a flow over a flat plate obstacle problem. The test is adopted from \cites{leotest,CIBIK2025117526}. The computational domain $\Omega$ consists of a rectangular channel $[-7,20]\times[0,20]$ with an embedded flat plate obstacle of size $0.125\times 1$. A domain with generated mesh shows in Figure \ref{fig:test3-mesh}. The inflow boundary condition is given by $\bfu_{\text{in}} = \langle 1,0 \rangle^\top$, generating a unidirectional flow in the positive $x$-direction. No-slip boundary conditions are imposed on the upper and lower walls as well as on the plate obstacle. A weak zero-traction condition is imposed at the outflow boundary as in \cite{leotest}.

We use the backward Euler scheme for the time discretization. The nonlinear convection term is treated using a skew-symmetric formulation, with the advecting velocity lagged at the previous time step, i.e. $b(\bfv^h_n,\bfv^h_{n+1}, \Theta)$. We employ an inf‑sup stable Taylor–Hood finite elements using $\mathbb{P}_2$ for velocity and $\mathbb{P}_1$ for pressure. The reference solution and the assimilated system are computed on a mesh with $56{,}790$ velocity degrees of freedom and $7{,}278$ pressure degrees of freedom, corresponding to a total of $64{,}068$ DOF. The reference solution is projected onto the $\mathbb{P}_I$ space via linear interpolation to obtain the coarse observation $I_H(\bfu)$. The data term is evaluated using lagged (previous time step) observations.
\subsection*{Simulation Setup} The assimilated solution is initialized from a quiescent flow, whereas the reference initial condition is obtained by solving the steady Stokes problem. We set $\nu=1$. The Reynolds number is defined as $Re = \frac{UL}{\nu}$, where the characteristic velocity $U$ is taken as the inlet velocity magnitude ($U = 1$), and the characteristic length $L$ is chosen as the length of the computational domain ($L = 27$). We set timestep size $\Delta t =0.02$ and the final time $T=7$. We investigate the effect of nonlinear nudging by comparing it with the classical linear nudging approach. The nonlinear method is tested with $\mu=5$, $\chi =1$  and $\gamma =0.1$, while the linear baseline is considered under two regimes: ($\mu=5$, $\chi =0$), and  ($\mu=5$, $\chi=0$). This allows us to distinguish the contribution of nonlinear feedback from that of increased linear nudging strength.

\subsection*{Results} 
\begin{figure}
    \centering
\includegraphics[width=0.7\linewidth]{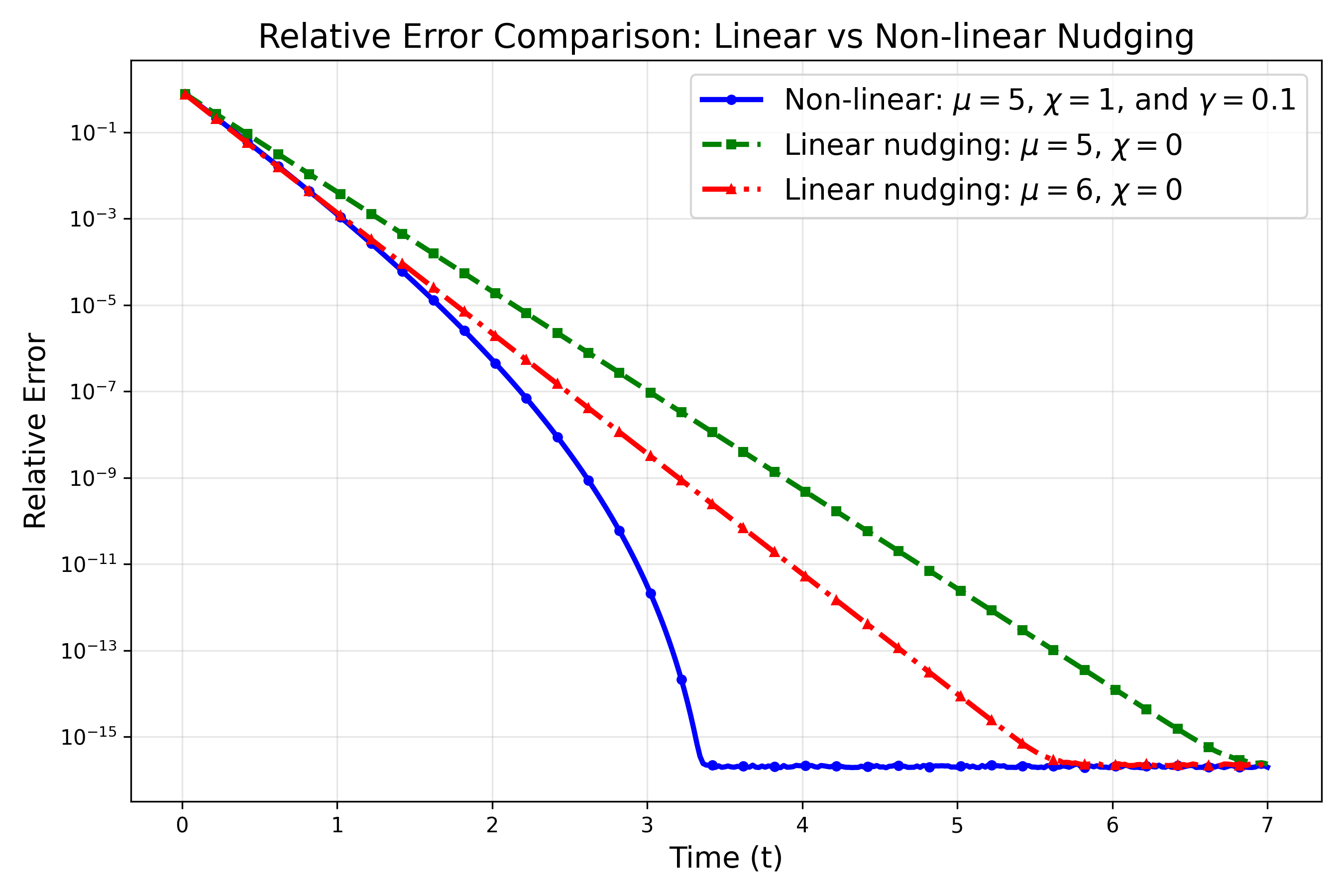}
\caption{Nonlinear nudging ($\mu=5$, $\chi=1$, $\gamma=0.1$) and linear nudging ($\mu=5,6$) all reach $\mathcal{O}(10^{-16})$ error. Initially, nonlinear nudging behaves like linear nudging with $\mu=6$, but after $t=1$ it decays faster, reaching steady state by $t\approx3.5$, while linear cases take longer.} 
    \label{fig:test3-comparison-non-linear}
\end{figure}

In Figure \ref{fig:test3-comparison-non-linear}, the nonlinear nudging case ($\mu = 5$, $\chi = 1$, $\gamma = 0.1$) and the linear nudging cases ($\mu = 5$ and $\mu = 6$) all reach a relative error of $\mathcal{O}(10^{-16})$ over a long time horizon. At early times, the nonlinear nudging exhibits behavior similar to linear nudging with $\mu = 6$. After $t = 1$, the error in the nonlinear nudging case decreases much more rapidly than in the linear cases, reaching a steady state around $t = 3.5$. In contrast, linear nudging with $\mu = 6$ requires a longer time (approximately $t = 6$) to reach the same steady state, while linear nudging with $\mu = 5$ takes even longer. 

\begin{figure}
        \centering
\includegraphics[width=0.7\linewidth]{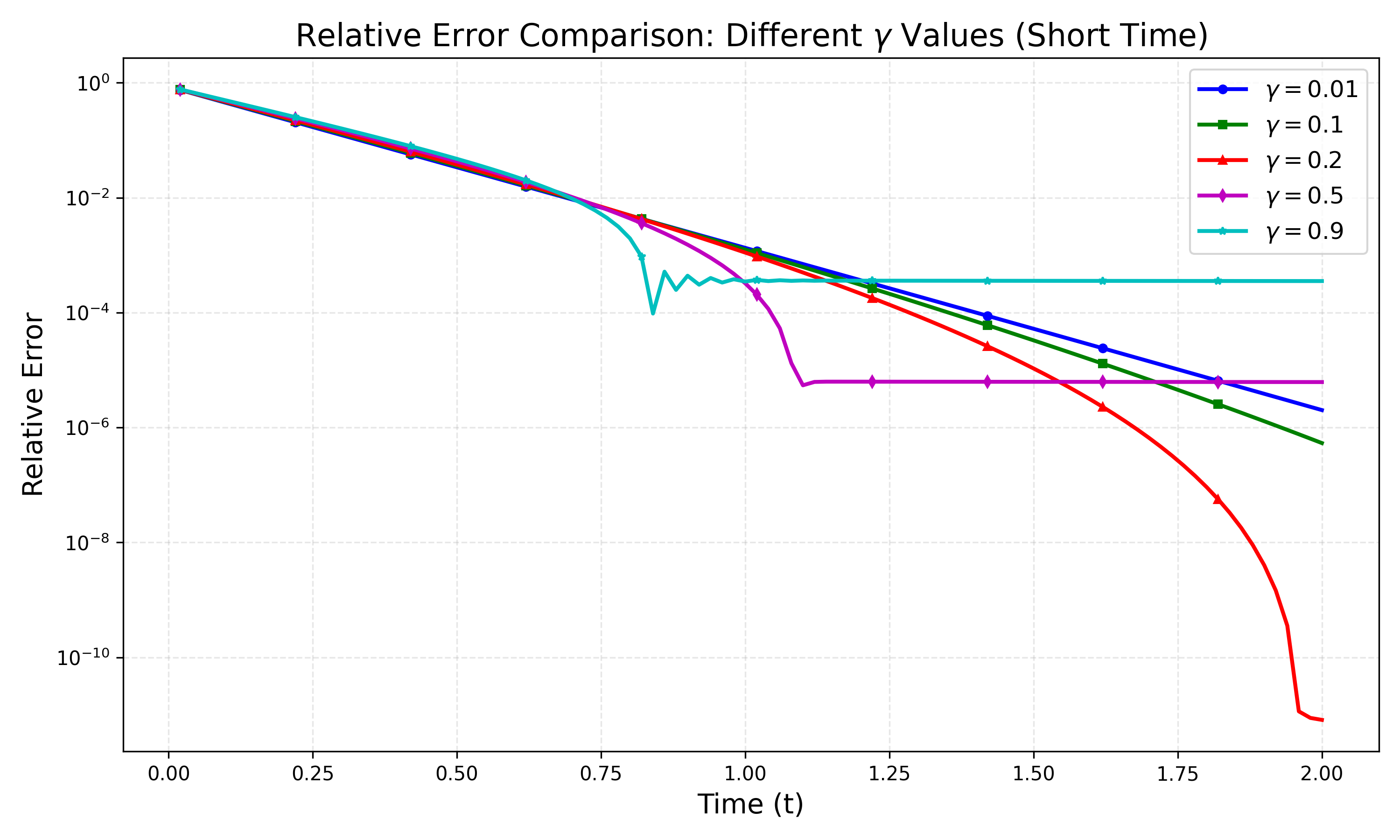}
        \caption{Relative velocity error for varying nonlinear nudging power $\gamma$, with $\mu = 5$ and $\beta = 1$, over a short time horizon $T = 2$. For small values of $\gamma$, the nonlinear nudging exhibits behavior similar to that of linear nudging with $\mu = 6$. As $\gamma$ increases, the error decays more rapidly compared to the linear case; it eventually reaches a plateau at a higher error level.}
        \label{fig:test3-different-gamma}
    \end{figure}

In Figure~\ref{fig:test3-different-gamma}, we study the effect of the nonlinear nudging power $\gamma$ over a short time horizon $T = 2$, with fixed parameters $\mu = 5$ and $\beta = 1$. For small values of $\gamma$, the nonlinear nudging exhibits behavior similar to linear nudging with $\mu = 6$. As $\gamma$ increases, the error decays more rapidly than in the linear case; however, it eventually reaches a higher plateau. Although no tolerance $\epsilon$ is prescribed in the numerical experiments, these observations are consistent with the sufficient mesh-resolution condition in equation~(\ref{H_for_nonlinear}). For $0<\gamma<1$, this condition can be written as
\begin{equation}
H \leq
\left(\frac{\alpha C_1^2 \mu\epsilon^2}{\nu}\right)^{\frac{\gamma}{2(1-\gamma)}}
\leq
\left(\frac{C_1^2\mu^2\epsilon^2}{\nu}\right)^{\frac{\gamma}{2(1-\gamma)}}
.
\end{equation}
In the small-tolerance regime considered in the analysis, assume that
$0<\frac{C_1^2\mu^2\epsilon^2}{\nu}<1$. Since $\gamma/[2(1-\gamma)]$ increases with $\gamma$, the quantity $\left(\frac{C_1^2\mu^2\epsilon^2}{\nu}\right)^{\frac{\gamma}{2(1-\gamma)}}$ decreases as $\gamma$ increases. Thus, larger values of $\gamma$ impose a more restrictive sufficient
upper bound on the mesh size $H$. Because $H$ is fixed in the numerical
experiments, this condition may cease to hold for sufficiently large
$\gamma$, providing a possible explanation for the observed higher
asymptotic error plateau.

\subsection{Flow past a cylinder}\label{sec:cylinder}
In this section, we test the nonlinear nudging method on the flow past a cylinder in a channel, a standard benchmark for incompressible flow solvers \cite{John2016}. Unlike the periodic box test, this problem has a physical domain with an obstacle, no-slip walls, an inflow and an outflow boundary, and its numerical scheme is the closest to Algorithm \ref{Algorithm} among our tests. We study two observation settings: data over the whole domain, and data only in a subregion containing the cylinder.

\subsection*{Simulation Setup}
\begin{figure}
    \centering
    \begin{subfigure}{0.95\textwidth}
        \centering
        \includegraphics[width=\linewidth]{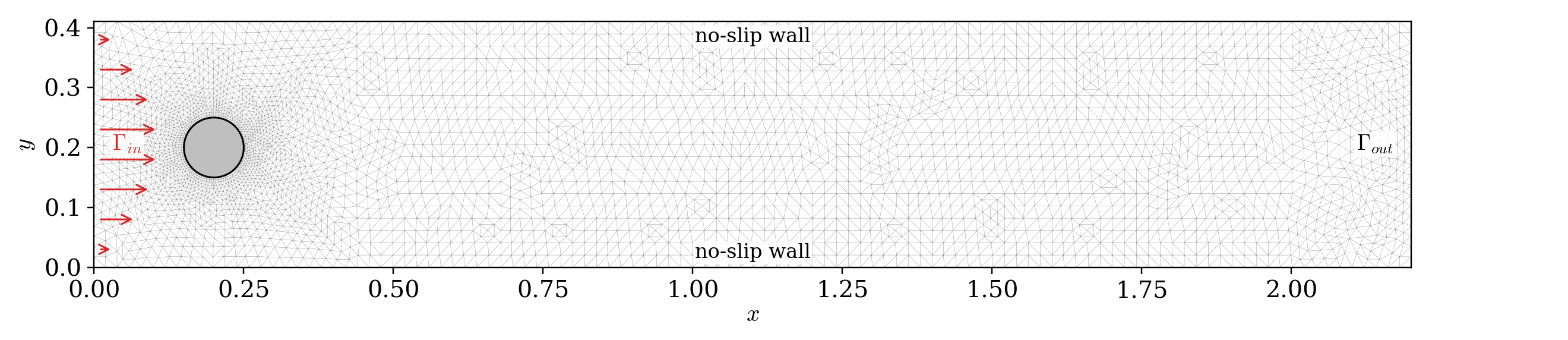}
        \caption{Domain and DNS mesh.}
        \label{fig:cyl_domain}
    \end{subfigure}
    \begin{subfigure}{0.95\textwidth}
        \centering
        \includegraphics[width=\linewidth]{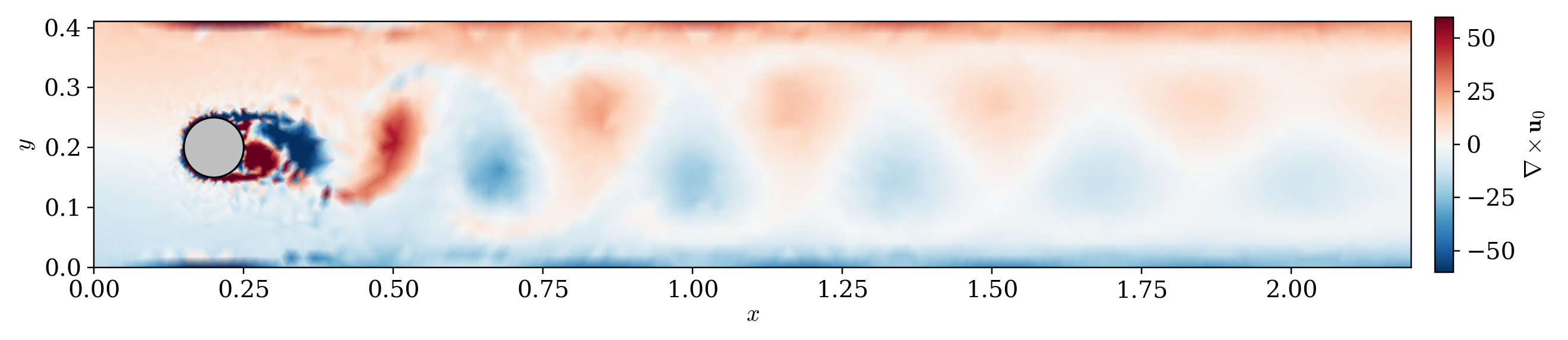}
        \caption{Vorticity of the initial condition $\bfu_0^h$.}
        \label{fig:cyl_u0}
    \end{subfigure}
    \caption{Flow past a cylinder, $Re=1000$.}
    \label{fig:cyl_setup}
\end{figure}

The domain is the channel $[0,2.2]\times[0,0.41]$ with a cylinder of diameter $D=0.1$ centered at $(0.2,0.2)$; see Figure \ref{fig:cyl_domain}. We impose the parabolic inflow $\bfu_{in}=\langle 4U_m y(0.41-y)/0.41^2,\,0\rangle^\top$ with $U_m=1.5$, i.e., a mean inflow velocity $\bar U=1$, no-slip conditions on the walls and on the cylinder, and the do-nothing condition at the outflow. We set $\nu=10^{-4}$, so that the Reynolds number based on $\bar U$ and $D$ is $Re=\bar U D/\nu=1000$.

The coarse observational mesh $\Omega_H$ is generated by Delaunay triangulation, and each coarse triangle is subdivided into four to obtain the fine DNS mesh $\Omega_h$, i.e., $H=2h$. We use Taylor--Hood $(P_2,P_1)$ elements, with $29{,}880$ velocity and $3{,}822$ pressure degrees of freedom. For time integration, we use the backward Euler scheme with the convecting velocity lagged, $(\bfv^h_n\cdot\nabla)\bfv^h_{n+1}$, and $\Delta t=0.005$, as in Algorithm \ref{Algorithm}; the denominator of the nonlinear nudging term is evaluated at the previous time level, as in the remark following Algorithm \ref{Algorithm}. The interpolant $I_H$ is the nodal piecewise linear interpolant on $\Omega_H$. 

To generate the initial condition $\bfu_0^h$ for the reference solution, we simulate the system from rest over $t\in[-6,0]$. The resulting flow is a time-periodic von Karman vortex street; see Figure \ref{fig:cyl_u0}. 
%with Strouhal number $St=fD/\bar U\approx0.34$, where $f$ is the shedding frequency; 

%To generate the initial condition $\bfu_0^h$ for the reference solution, we simulate the system from rest over $t\in[-6,0]$. The resulting flow is a periodic von K\'arm\'an vortex street; see Figure \ref{fig:cyl_u0}. The drag and lift coefficients, $C_D=2F_x/(\bar U^2D)$ and $C_L=2F_y/(\bar U^2D)$, where $(F_x,F_y)$ is the force exerted by the fluid on the cylinder, are periodic in time with mean drag $\bar C_D=2.586$, maximum drag $C_{D,\max}=3.074$, maximum lift $|C_L|_{\max}=2.725$, and Strouhal number $St=fD/\bar U=0.344$, where $f$ is the shedding frequency; these values are computed over $t\in[40,50]$.

The assimilated solution is initialized with $\bfv^h_0=\mathbf{0}$, and all simulations are run up to $T=50$. We compare the following cases:
\begin{itemize}
    \item linear nudging with $\mu=1$ ($\chi=0$),
    \item nonlinear nudging with $\mu=\chi=1$ and $\gamma=0.1,\ 0.25,\ 0.5,\ 0.9$.
\end{itemize}

\subsection*{Results: observations over the whole domain}
\begin{figure}
    \centering
    \includegraphics[width=0.7\linewidth]{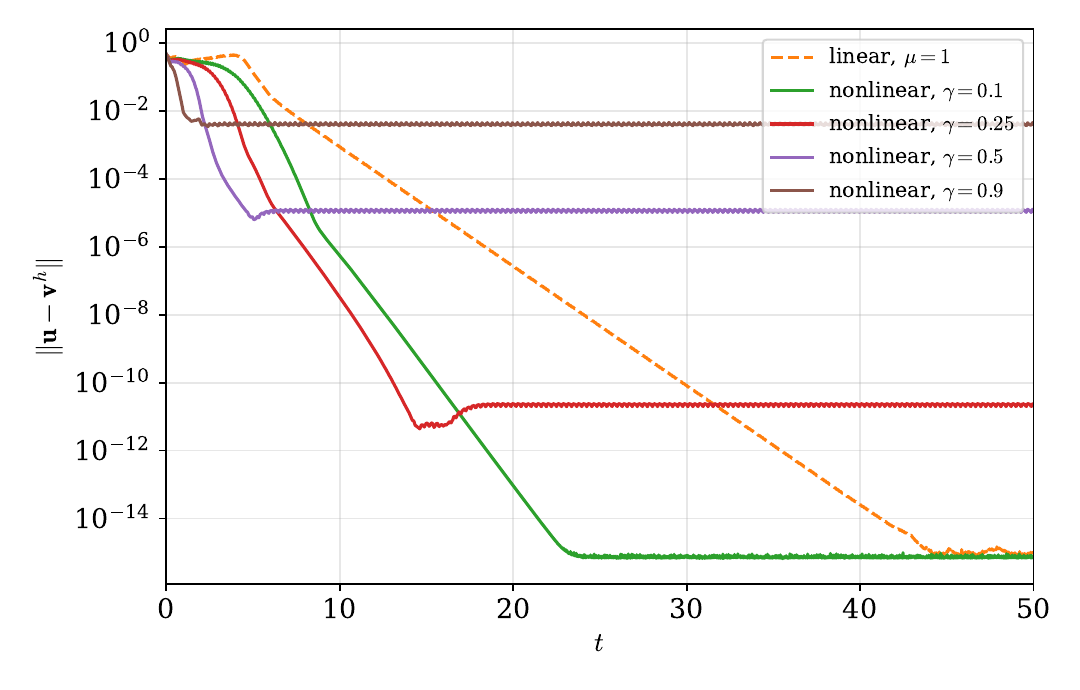}
    \caption{Flow past a cylinder, observations over the whole domain: error $\|\bfu-\bfv^h\|$ for linear nudging ($\mu=1$, $\chi=0$) and nonlinear nudging ($\mu=\chi=1$).}
    \label{fig:cyl_global}
\end{figure}

Figure \ref{fig:cyl_global} shows the error $\|\bfu_n-\bfv^h_n\|$ when the data are available over the whole domain. Linear nudging decays exponentially at a rate of about $0.8$ and reaches machine precision at about $t=45$. The nonlinear case $\gamma=0.1$ converges about twice as fast: the error drops below $10^{-8}$ at $t=13$, compared with $t=25$ for the linear case, and reaches machine precision at about $t=23$. Larger values of $\gamma$ give a faster decay at early times (the error drops below $10^{-4}$ at $t=5.5$ and $t=3.5$ for $\gamma=0.25$ and $\gamma=0.5$), but the error then stagnates at a plateau: $2.5\times10^{-11}$, $10^{-5}$, and $4.0\times10^{-3}$ for $\gamma=0.25$, $0.5$, and $0.9$, respectively.  In all three cases, the plateau is reached when the effective coefficient $\mu+\chi\|I_H(\bfu_n-\bfv^h_n)\|^{-\gamma}$ attains approximately $550$--$700$; denoting this value by $\mu_*$, the corresponding observed error, $\|I_H(\bfu_n-\bfv^h_n)\|\approx\big((\mu_*-\mu)/\chi\big)^{-1/\gamma}$, increases with $\gamma$, which explains the ordering of the plateaus. We do not observe such a limit in the periodic box test, where the effective coefficient reaches $2.3\times 10^{7}$, and we leave its explanation to future work.

\subsection*{Results: local observations}
\begin{figure}
    \centering
    \begin{subfigure}{0.95\textwidth}
        \centering
        \includegraphics[width=\linewidth]{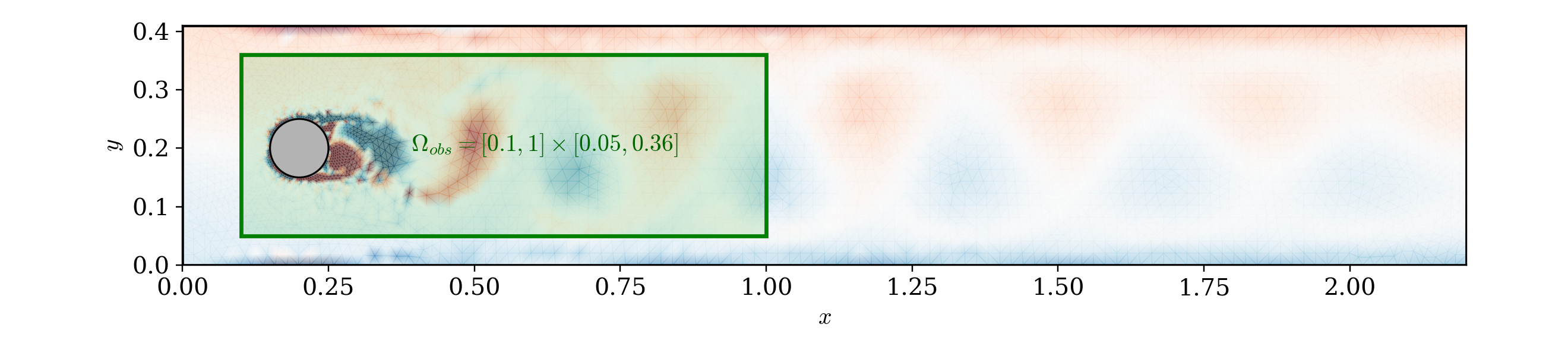}
        \caption{Observational region $\Omega_{obs}$ (shaded) which  covers $ \sim 30 \%$ of the domain, shown over the vorticity of $\bfu_0^h$.}
        \label{fig:cyl_local_region}
    \end{subfigure}
    \begin{subfigure}{0.7\textwidth}
        \centering
        \includegraphics[width=\linewidth]{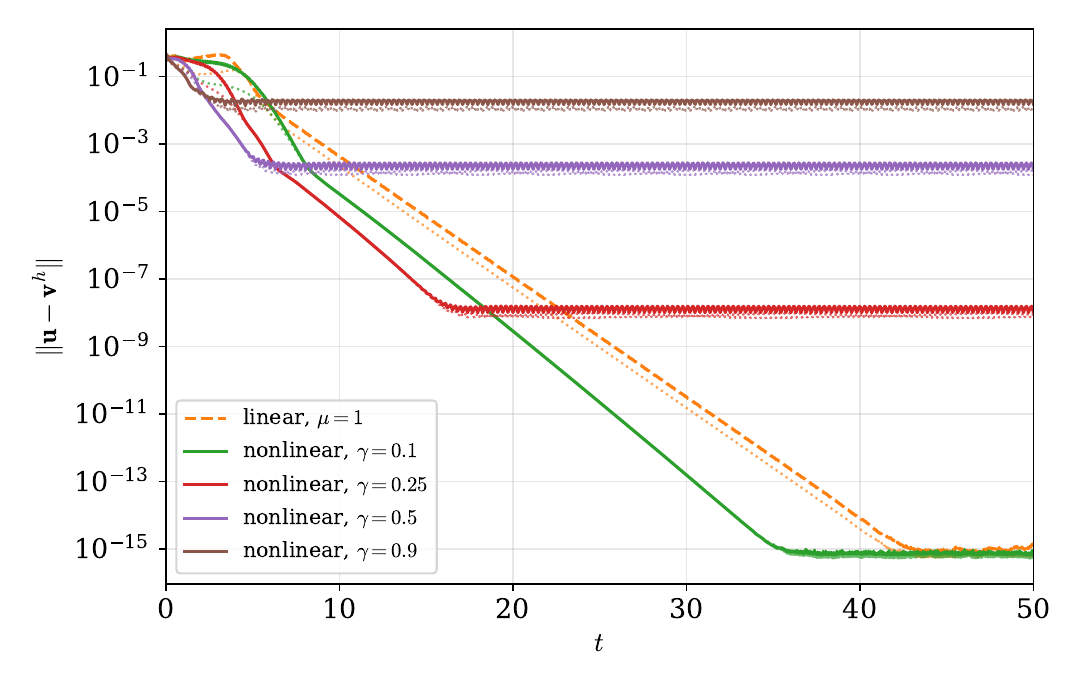}
        \caption{Error $\|\bfu-\bfv^h\|$ over $\Omega$ (solid) and over $\Omega\setminus\Omega_{obs}$ (dotted).}
        \label{fig:cyl_local_err}
    \end{subfigure}
    \caption{Flow past a cylinder, local observations in $\Omega_{obs}$: linear nudging ($\mu=1$, $\chi=0$) and nonlinear nudging ($\mu=\chi=1$).}
    \label{fig:cyl_local}
\end{figure}

Motivated by the recent work \cite{FangPakzad2026},  we next assume that observations are available only in the subregion $\Omega_{obs}=[0.1,1]\times[0.05,0.36]$, which contains the cylinder and the near wake and covers about $30\%$ of the domain; see Figure \ref{fig:cyl_local_region}. The nudging terms are multiplied by the indicator function of $\Omega_{obs}$, and the norm in the denominator of the nonlinear term is taken over $\Omega_{obs}$. The error is measured over the whole domain $\Omega$.

Figure \ref{fig:cyl_local_err} shows that the local data are sufficient to synchronize the whole flow, including the far wake downstream of $\Omega_{obs}$: the error outside $\Omega_{obs}$ (dotted lines) follows the total error. The behavior is similar to the case with observations over the whole domain. Linear nudging and nonlinear nudging with $\gamma=0.1$ reach machine precision at about $t=44$ and $t=35$, respectively, while $\gamma=0.25$, $0.5$, and $0.9$ reach plateaus of $1.0\times 10^{-8}$, $1.8\times10^{-4}$, and $1.5\times10^{-2}$. Here, the plateaus are reached when the effective coefficient attains approximately $140$. 

The location of the observational region matters. When $\Omega_{obs}=[0.3,1]\times[0.05,0.36]$ (the test is not included here),  which excludes the cylinder, linear nudging still converges, to $1.8\times10^{-6}$ at $t=50$, but the nonlinear cases $\gamma=0.1$ and $0.25$ do not synchronize.  

\subsection*{Results: Drag and lift}

We also compare the drag and lift coefficients, $C_D$ and $C_L$, of the assimilated solutions with those of the reference solution. These coefficients are the horizontal and vertical components of the force exerted by the fluid on the cylinder, normalized by $\bar U^2D/2$, and are computed in the same way for both solutions,  see \cite{SchaferTurek1996} for details. Unlike the error $\|\bfu-\bfv^h\|$, which is a global quantity, they depend on the velocity gradient and the pressure at the cylinder surface, and therefore measure the accuracy of the assimilated solution near the obstacle.

Table \ref{tab:cyl_forces} reports the maximum difference between the coefficients of the assimilated and the reference solutions over $t\in[40,50]$, compared at the same time instants. For the cases that reach machine precision, the coefficients agree up to the output precision, $10^{-10}$. For the cases with a plateau, the differences increase with $\gamma$, up to $1.3\times10^{-2}$ in $C_L$ for $\gamma=0.9$ with local observations. All assimilated solutions shed vortices at the same frequency as the reference solution. In contrast, when $\Omega_{obs}=[0.3,1]\times[0.05,0.36]$, the cases $\gamma=0.1$ and $0.25$ differ by up to $3.8$ in $C_L$, which exceeds the maximum lift of the reference solution, $2.7$; their lift is out of phase with that of the reference solution.

\begin{table}[ht]
    \centering
    \setlength{\tabcolsep}{8pt}
    \begin{tabular}{|c|c|c|c|c|}
        \hline
        & \multicolumn{2}{c|}{whole domain} & \multicolumn{2}{c|}{local, $\Omega_{obs}$} \\
        \hline
        Case & $\Delta C_D$ & $\Delta C_L$ & $\Delta C_D$ & $\Delta C_L$ \\
        \hline
        linear, $\mu=1$          & $<10^{-10}$ & $<10^{-10}$ & $<10^{-10}$ & $<10^{-10}$ \\
        nonlinear, $\gamma=0.1$  & $<10^{-10}$ & $<10^{-10}$ & $<10^{-10}$ & $<10^{-10}$ \\
        nonlinear, $\gamma=0.25$ & $<10^{-10}$ & $<10^{-10}$ & $3.9\times10^{-9}$ & $8.6\times10^{-9}$ \\
        nonlinear, $\gamma=0.5$  & $7.5\times10^{-7}$ & $9.7\times10^{-6}$ & $7.6\times10^{-5}$ & $1.5\times10^{-4}$ \\
        nonlinear, $\gamma=0.9$  & $2.5\times10^{-4}$ & $3.7\times10^{-3}$ & $7.5\times10^{-3}$ & $1.3\times10^{-2}$ \\
        \hline
    \end{tabular}
    \caption{Flow past a cylinder: maximum deviation of the drag and lift coefficients of the assimilated solutions from those of the reference solution over $t\in[40,50]$, with observations over the whole domain and in $\Omega_{obs}=[0.1,1]\times[0.05,0.36]$. Nonlinear cases use $\mu=\chi=1$.}
    \label{tab:cyl_forces}
\end{table}

\section{Discussion}

In this work, we developed and analyzed a fully discrete scheme for the nonlinear data assimilation system applied to the incompressible Navier--Stokes equations \eqref{NLDANSE}. We established unconditional long-time
stability for the fully discrete scheme and derived error estimates for the continuous in time and discrete in space scheme. The finite element method is employed for the spatial discretization, allowing flexibility in handling complex geometries and boundary conditions. For computational efficiency and ease of implementation, we consider the first-order time discretization. Establishing error estimates for the fully discrete scheme and
developing and analyzing higher-order time discretizations remain
open problems within the present framework.

An important avenue for future work is the extension of the present framework to three-dimensional simulations. In addition to the increased computational complexity 
of three-dimensional flows, such an extension would, from a theoretical standpoint, require assuming the existence of sufficiently regular solutions, along with addressing 
additional stability and analytical challenges.

\appendix
\section{Stability analysis for a general interpolation operator (\ref{Ineq:Intrpolate})}
\label{app:general-interpolant}
\begin{thm}\label{stability-g}
Let $0\le \gamma<1$, and assume that
$U:=\|\bfu\|_{L^{\infty}(0,\infty;L^2(\Omega))}<\infty$ and
$F:=\|\bff\|^2_{L^{\infty}(0,\infty;H^{-1}(\Omega))}<\infty$. If
$$\mu C_1^2H^2\leq \frac{\nu}{4},  \quad\text{and}\quad \chi\leq\mu,$$
then Algorithm \ref{Algorithm} is stable, i.e.,  for all $n\geq 0$,
\begin{equation}\label{stability-final-result-g}
\|\bfv^h_{n+1}\|^2 \leq
\Big(1+\Delta t\Big(\frac{\nu}{C_{FP}^2}+\frac{\mu}{2}\Big)\Big)^{-(n+1)}\|\bfv^h_0\|^2
+\frac{M}{\frac{\nu}{C_{FP}^2}+\frac{\mu}{2}},
\end{equation}
where
$M:=\frac{2F}{\nu}+2\mu C_I^2U^2+4\chi C_I^{2-2\gamma}U^{2-2\gamma}
+\chi\gamma\big(4(2-\gamma)\big)^{(2-\gamma)/\gamma}C_I^{\frac{2-2\gamma}{\gamma}}$
is independent of $\Delta t$, $h$, and $n$.
\end{thm}

\begin{proof}
We consider the case $0<\gamma<1$. When $\gamma=0$,
Algorithm~\ref{Algorithm} reduces to the standard linear-nudging
scheme. We set $\Theta^{h} = \mathbf{v}^h_{n+1}$ in \ref{Algorithm}. We add and subtract $\mu\, \|v_{n+1}^h\|^2$ within the initial component of the linear nudging term. Through the application of the polarization identity, and by the skew-symmetry of the nonlinear term, we reach the following expression:
\begin{equation}\label{Disc1-A}
\begin{gathered}
\frac{1}{2\Delta t } \left( \| \bfv^h_{n+1}\|^2 -\| \bfv^h_{n}\|^2  \right)  + \nu\, \|\nabla \bfv^h_{n+1}\|^2  + \mu \,  \|\bfv^h_{n+1}\|^2 \\\leq (\bff_{n+1} , \bfv^h_{n+1}) 
 - \mu \, (I_H(\bfv^h_{n+1})- \bfv^h_{n+1},  \bfv^h_{n+1}) \\
 + \mu \, (I_H(\bfu_{n+1}),  \bfv^h_{n+1})  + \chi \,  \left( \frac{I_H(\bfu_{n+1}) -  I_H(\bfv^h_{n+1})}{\|(I_H(\bfu_{n+1}) -  I_H(\bfv^h_{n+1}))\|^{\gamma}},  \bfv^h_{n+1}\right).
\end{gathered}
\end{equation}

Similarly as in Theorem \ref{stability}
$$(\bff_{n+1} , \bfv^h_{n+1}) \leq \frac{1}{\nu}\|\bff_{n+1}\|^2_{-1} + \frac{\nu}{4} \| \nabla \bfv^h_{n+1}\|^2.$$
Applying  Cauchy-Schwarz inequality, Young's inequality, and interpolation property \eqref{Ineq:Intrpolate}  to the second term,  along with the assumption  $\mu C_1^2 H^2 \leq \nu/4$,  resulting in
\begin{equation*}
\begin{split}
\mu \, (I_H(\bfv^h_{n+1})- \bfv^h_{n+1},  \bfv^h_{n+1})
&\leq \mu \|I_H\bfv^h_{n+1} - \bfv^h_{n+1}\| \|\bfv^h_{n+1}\|\\
&\leq \mu\,C_1^2\,  H^2 \|\nabla \bfv_{n+1}^h\|^2 + \frac{\mu}{4} \|\bfv_{n+1}^h\|^2\\
& \leq   \frac{\nu}{4} \|\nabla \bfv_{n+1}^h\|^2 + \frac{\mu}{4} \|\bfv_{n+1}^h\|^2. 
\end{split}
\end{equation*}
Using \eqref{Ineq:Intrpolate},  one can estimate the third term as 
$$\mu \, (I_H(\bfu_{n+1}),  \bfv^h_{n+1}) \leq \mu \, C_I^2 \|\bfu_{n+1}\|^2 + \frac{\mu}{4}\, \|\bfv_{n+1}^h\|^2.$$
We ultimately employ Cauchy-Schwarz inequality and Young's inequality with the interpolation stability bound \eqref{Ineq:Intrpolate}. Now we bound the nonlinear nudging term in equation \eqref{Disc1-A} as
\begin{equation}
\begin{gathered}
\chi \,  \left( \frac{I_H(\bfu_{n+1}) -  I_H(\bfv^h_{n+1})}{\|I_H(\bfu_{n+1}) -  I_H(\bfv^h_{n+1})\|^{\gamma}},  \bfv^h_{n+1}\right)   \leq \chi  \|I_H(\bfu_{n+1}) -  I_H(\bfv^h_{n+1})\|^{1-\gamma}\, \| \bfv^h_{n+1}\|\\
\leq \chi \|I_H(\bfu_{n+1}) \|^{1-\gamma}\, \| \bfv^h_{n+1}\| + \chi \|I_H(\bfv^h_{n+1}) \|^{1-\gamma}\, \| \bfv^h_{n+1}\|\\
 \leq \chi\,  C_I^{1-\gamma}\, U^{1-\gamma}\, \|\bfv^h_{n+1}\| + \chi\,  C_I^{1-\gamma}\,  \|\bfv^h_{n+1}\|^{2-\gamma}\\
 \leq 2\chi C_I^{2-2\, \gamma}\, U^{2-2\, \gamma} + \frac{\chi}{8}  \|\bfv^h_{n+1}\|^2  + \chi \frac{\gamma}{2} \left(4(2-\gamma)\right)^{(2-\gamma)/\gamma} C_I^{\frac{2 - 2 \gamma}{\gamma}} + \frac{\chi}{8} \|\bfv^h_{n+1}\|^2 \\
 \leq   2\chi C_I^{2-2\, \gamma}\, U^{2-2\, \gamma}  +  \chi \frac{\gamma}{2} \left(4(2-\gamma)\right)^{(2-\gamma)/\gamma} C_I^{\frac{2 - 2 \gamma}{\gamma}} + \frac{\chi}{4} \| \bfv^h_{n+1}\|^2.
\end{gathered}
\end{equation}
Combining the above bounds for the right hand side of \eqref{Disc1-A} and multiplying by $2\Delta t$, we have
\begin{equation*}
\begin{split}
\left( \| \bfv^h_{n+1}\|^2 -\| \bfv^h_{n}\|^2  \right)  + \Delta t \, \nu\, \|\nabla \bfv^h_{n+1}\|^2  + \Delta t \, \frac{1}{2}\mu \,  \|\bfv^h_{n+1}\|^2 +\Delta t \frac{1}{2}(\mu-\chi) \|\bfv^h_{n+1}\|^2\\
\leq  \Delta t \left( \frac{2F}{\nu} + 2\mu \, C_I^2 U^2 +  4 \chi  C_I^{2-2\, \gamma}\, U^{2-2\, \gamma}  +  \chi \gamma \left(4(2-\gamma)\right)^{(2-\gamma)/\gamma} C_I^{\frac{2 - 2 \gamma}{\gamma}} \right).
\end{split}
\end{equation*}
Since $\chi \leq \mu$, thus
\begin{equation}
\begin{gathered}
\left( \| \bfv^h_{n+1}\|^2 -\| \bfv^h_{n}\|^2  \right)  + \Delta t \, \nu\, \|\nabla \bfv^h_{n+1}\|^2  + \Delta t \, \frac{1}{2}\mu \,  \|\bfv^h_{n+1}\|^2 \\
\leq  \Delta t \left( \frac{2F}{\nu} + 2\mu \, C_I^2 U^2 +  4 \chi  C_I^{2-2\, \gamma}\, U^{2-2\, \gamma}  +  \chi \gamma \left(4(2-\gamma)\right)^{(2-\gamma)/\gamma} C_I^{\frac{2 - 2 \gamma}{\gamma}} \right).
\end{gathered}
\end{equation}
Applying the Poincar\'e inequality to the viscous term and collecting terms gives us
\begin{equation}
\begin{gathered}
\left( 1 + \Delta t \, \nu \frac{1}{(C_{PF})^2} + \Delta t\,  \frac{1}{2}\mu \right)  \| \bfv^h_{n+1}\|^2  \leq   \| \bfv^h_{n}\|^2  \\
+\Delta t \left( \frac{2F}{\nu} + 2 \mu \, C_I^2 U^2 +  4 \chi  C_I^{2-2\, \gamma}\, U^{2-2\, \gamma}  +  \chi \gamma \left(4(2-\gamma)\right)^{(2-\gamma)/\gamma} C_I^{\frac{2 - 2 \gamma}{\gamma}} \right). 
\end{gathered}
\end{equation}
Using Lemma \ref{Lem;Seq}, we have
\begin{equation}
\begin{gathered}
\| \bfv^h_{n+1}\|^2    \leq   \frac{1}{\left( 1 + \Delta t \, \nu \frac{1}{(C_{PF})^2} + \Delta t\,  \frac{1}{2}\mu \right)^{1+n}}\, \| \bfv^h_{0}\|^2\\
+ \frac{  \frac{2F}{\nu} + 2\mu \, C_I^2 U^2 +  4 \chi  C_I^{2-2\, \gamma}\, U^{2-2\, \gamma}  +  \chi \gamma \left(4(2-\gamma)\right)^{(2-\gamma)/\gamma} C_I^{\frac{2 - 2 \gamma}{\gamma}} }{ \Delta t \, \nu \frac{1}{(C_{PF})^2} + \Delta t\,  \frac{1}{2}\mu} \leq C (\text{data}).
 \end{gathered}
\end{equation}
\end{proof}
\section*{Acknowledgment}
A.~P. was  supported by NSF grant DMS-2532987.  A.~P.  also acknowledges the
Indiana University Pervasive Technology Institute for providing HPC (Big Red 3, Carbonate) and storage resources that have contributed to the research results reported in this article.

\bibliographystyle{plain}
\bibliography{mybib}

\end{document}